\documentclass{article}
\usepackage{amsmath}
\usepackage{amsfonts}
\usepackage{subcaption}
\usepackage{graphicx}
\usepackage{float}
\usepackage[dutch,english]{babel}
\usepackage{amsmath,amssymb}
\usepackage{amsthm}
\usepackage[usenames,dvipsnames]{color}
\usepackage{accents}
\usepackage{xcolor}
\usepackage{mathtools}
\usepackage{multicol}
\usepackage{comment}
\usepackage{caption}
\usepackage[inline]{enumitem}
\usepackage[a4paper, left=35mm, right=35mm, top=35mm, bottom=35mm]{geometry}

\usepackage[
    backend=biber,
    style=numeric,
    sorting=none,
    url=false
]{biblatex}
\usepackage{csquotes}
\renewbibmacro{in:}{%
	\ifentrytype{article}{}{\printtext{\bibstring{in}\intitlepunct}}}

\allowdisplaybreaks

\usepackage{standalone}
\usepackage{tikz}
\usetikzlibrary{shadows}
\usetikzlibrary{arrows}
\usetikzlibrary{arrows.spaced}
\usetikzlibrary{arrows.meta}
\usetikzlibrary{shapes.geometric}
\usetikzlibrary{fit}
\usetikzlibrary{backgrounds}
\usetikzlibrary{calc}

\definecolor{red}{RGB}{215,51,103}
\definecolor{blue}{RGB}{24,176,226}
\definecolor{grey}{RGB}{85,85,85}
\definecolor{mix}{RGB}{120,114,165}
\definecolor{uhhdarkgrey}{RGB}{142,169,214}
\definecolor{uhhsteingrau}{RGB}{59,81,91}
\definecolor{weirdgrey}{RGB}{105,82,64}
\definecolor{notevengrey}{RGB}{92,97,21}
\definecolor{upborange}{RGB}{242,149,18}
\definecolor{newmix}{RGB}{160,125,116}
\definecolor{mycolor}{RGB}{190, 0, 190}
\definecolor{mycolor2}{RGB}{190, 80, 100}
\definecolor{wred}{rgb}{0.7,0.18,0.12}
\definecolor{wgreen}{rgb}{0.1,0.53,0.37}

\usepackage[pdftex,hidelinks]{hyperref}
\usepackage[capitalise,nameinlink]{cleveref}

\newcommand{\la}{\lambda}
\newcommand{\NN}{\mathbf{N}}
\newcommand{\NNd}{\NN^{\diamondsuit}}

\theoremstyle{plain}
\newtheorem{thr}{Theorem}[section]
\Crefname{thr}{Theorem}{Theorems}

\newtheorem{lem}[thr]{Lemma}
\Crefname{lem}{Lemma}{Lemmas}

\Crefname{cor}{Corollary}{Corollaries}

\theoremstyle{definition}
\newtheorem{defi}[thr]{Definition}
\crefname{defi}{definition}{definitions}
\newtheorem{ex}[thr]{Example}
\Crefname{ex}{Example}{Examples}
\theoremstyle{remark}
\newtheorem{remk}{Remark}
\theoremstyle{remark}

\newcommand{\field}[1]{\mathbb{#1}}
\newcommand{\R}{\field{R}}
\newcommand{\N}{\field{N}}

\DeclareMathOperator{\sgn}{sgn}

\numberwithin{equation}{section}

\def\barray{\begin{eqnarray*}}             \def\earray{\end{eqnarray*}}
\def\beq{\begin{equation}} \def\eeq{\end{equation}}

\DeclarePairedDelimiter\ceil{\lceil}{\rceil}

\title{Higher order interactions lead to ``reluctant'' synchrony breaking}
\date{}
\author{S\"oren von der Gracht\thanks{Institute of Mathematics, Paderborn University, Germany, \href{mailto:soeren.von.der.gracht@uni-paderborn.de}{soeren.von.der.gracht@uni-paderborn.de}}, Eddie Nijholt
\thanks{\mbox{Department of Mathematics, Imperial College London, United Kingdom, \break \href{mailto:eddie.nijholt@gmail.com}{eddie.nijholt@gmail.com}}}, Bob Rink\thanks{\mbox{Department of Mathematics, Vrije Universiteit Amsterdam, The Netherlands, \href{mailto:b.w.rink@vu.nl}{b.w.rink@vu.nl}}}}

\excludecomment{percent}

\begin{document}

\maketitle

\begin{abstract}
To model dynamical systems on networks with higher order (non-pairwise) interactions, we recently introduced a new class of ODEs on hypernetworks. Here we consider one-parameter synchrony breaking bifurcations in such ODEs. We call a synchrony breaking  steady state branch ``reluctant'' if it is tangent to a synchrony space, but does not lie inside it.
We prove that reluctant synchrony breaking is ubiquitous in hypernetwork systems, by constructing a large class of examples that support it. We also give an explicit formula for the order of tangency to the synchrony space of a reluctant steady state branch.
\end{abstract}

\section{Introduction}
Recent advances in a large variety of research fields have highlighted the importance of non-pairwise interactions for the collective dynamical behavior of complex network systems. These so-called \emph{higher order interactions} turn out to be crucial in problems from, e.g.,  neuroscience (see \cite{Ariav.2003,Parastesh.2022}), social science (see \cite{Neuhauser.2022}), and ecology (see \cite{Chatterjee.2022,Gibbs.2022,Pal.2023}). Physicists have, in particular, emphasized the impact of group interactions on synchronisation behaviour in various coupled oscillator models and their generalizations (such as topological signals on cell complexes or swarmalators with higher order interactions) by comparing them to ``classical'' dyadic network models (e.g.~\cite{Carletti.2023,Anwar.2024}). At the same time, model reconstruction of dyadic networks is sometimes known to paradoxically yield group interactions, e.g. in experiments with electrochemical oscillators \cite{Nijholt.2022b}. Higher order interaction networks have consequently found their way into various recent mathematical studies as well. We mention in particular the theoretical papers  \cite{Mulas.2020,DeVille.2021b,Salova.2021b,Salova.2021,Nijholt.2022c,Aguiar.2022,vonderGracht.2023}, 
which investigate synchronisation in classes of networks with non-pairwise, nonlinear interactions in their equations of motion.
We also refer to the excellent surveys \cite{Battiston.2020,Porter.2020,Torres.2021,Bick.2023c,Boccaletti.2023} and references therein, for an in-depth discussion of higher order networks, and numerous examples of higher order network systems  arising in applications.

This paper builds on previous work of the authors   \cite{vonderGracht.2023}, in which we generalised the notion of a {\it coupled cell network},   introduced by Golubitsky, Stewart, Field et al. \cite{Golubitsky.2004, Golubitsky.2006, Field.2004}, to the context of higher order networks.  
We did this by introducing a class of ``hypernetworks'' and defining their ``admissible'' maps and ODEs, thus  
formalising the notion of a dynamical system on a higher order interaction network. 
We also introduced balanced colorings \cite{Golubitsky.2004} of hypernetworks, and hypergraph fibrations \cite{Boldi.2002, DeVille.2015}, and  used these concepts to classify the {\it robust synchrony patterns}, that is, the synchrony spaces that are invariant under every admissible map, to  hypernetwork dynamical systems.

The most surprising result in \cite{vonderGracht.2023} is the observation that the robust synchrony spaces of a hypernetwork system are not determined by linear terms in its equations of motion. This distinguishes hypernetworks from classical (dyadic) coupled cell networks, for which it was proved in \cite{Golubitsky.2005} that a synchrony space is  invariant under every admissible map, if and only if it is invariant under every linear admissible map, see also \cite{Aguiar.2014}. 
On the contrary, we prove in \cite{vonderGracht.2023} that a synchrony space of a hypernetwork system is robustly invariant, whenever it is invariant under all  polynomial admissible maps of a specific degree, which depends on the order of the hyperedges in the hypernetwork.
Examples moreover show that our estimate for this polynomial degree is sharp.

As a consequence, a hypernetwork-admissible map of sufficiently low polynomial degree may admit ``ghost'' synchrony spaces that are not supported by general, e.g., higher degree polynomial admissible maps. These ghost synchrony spaces may have a profound effect on the dynamics of the hypernetwork system. In particular, the final section of \cite{vonderGracht.2023} presents numerical evidence that they can give rise to a remarkable new type of local synchrony breaking bifurcation.  The aim of this paper is to explain when and why such bifurcations occur.

We in fact observed this type of bifurcation in a one-parameter family of admissible ODEs for the hypernetwork depicted in \Cref{fig:1st_ex}, meaning that these ODEs are of the form given in \Cref{systembifwithpara} below. 
\Cref{firstnumintro} displays two numerically obtained branches of steady states that emerge in  a bifurcation in a particular system of this form. The steady state branches were found by forward integrating the equations of motion -- so they are asymptotically stable. We see that $y_0=y_1$ for negative values of the bifurcation parameter $\lambda$, so on the negative branch $y_0$ and $y_1$ are synchronous. On the positive branch, $y_0$ and $y_1$ are nonsynchronous, i.e., for positive values of $\lambda$ it holds that $y_0\neq y_1$. However, even though $y_0$ and $y_1$ grow notably as   $\lambda$ increases, the difference $y_0-y_1$ only increases very slowly as a function of $\lambda$, and hence it appears that the branch is tangent to the  synchrony space $\{y_0=y_1\}$.  
In \cite{vonderGracht.2023}, we called this phenomenon ``reluctant synchrony breaking''. The term ``reluctant'' is to be understood literally, and refers to the slow separation of the states of two nodes that were synchronous before the bifurcation. Importantly, the reluctance is not caused by any (external) physical influence. On the contrary, it is solely due to the topology of the interaction structure of the hypernetwork, as we show below. A more detailed numerical analysis, see \Cref{firstnum-b,firstnum-c}, suggests that $y_0-y_1 \sim \lambda^3$, i.e., that the branch has a third order tangency to the synchrony space. In \Cref{sec:reluctant}, we  prove that this is indeed the case. 

We show in this paper that reluctant synchrony breaking is ubiquitous in hypernetworks. The main result that we present is \Cref{mainbifurcation}, which states that reluctant synchrony breaking occurs generically in one-parameter bifurcations in a large class of hypernetworks. These so-called {\it augmented hypernetworks} are constructed by coupling new nodes to an existing network or hypernetwork by means of specific higher order interactions. 
The hypernetwork depicted in \Cref{fig:1st_ex} is just one example of such an augmented hypernetwork. This means that the anomalous bifurcation that was discovered in \cite{vonderGracht.2023} and described above is not a numerical artefact. Instead, reluctant synchrony breaking is a generic phenomenon in ODEs of the form \eqref{systembifwithpara}. To illustrate our main result, we present several more examples in this paper. We also argue (see   {\rm \Cref{remk:secondclosing}}) that one may design augmented hypernetworks which admit reluctant synchrony breaking bifurcation branches with an arbitrarily high order of reluctancy, i.e., an arbitrarily high order of tangency to a synchrony space.

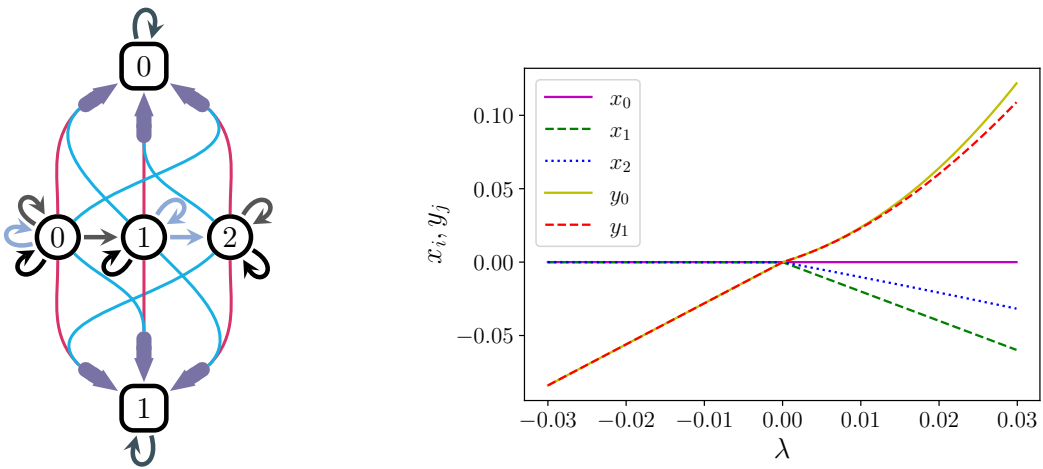
\begin{figure}
	\centering
	\begin{subfigure}{.35\linewidth}
		\resizebox{.8\linewidth}{!}{
			\begin{tikzpicture}[
	square/.style = {
		regular polygon,
		regular polygon sides=4
	},
	main node/.style={
		line width=1.5pt, 
		circle,
		draw,
		font=\sffamily,
		inner sep=2pt,
		fill=white
	},
	second node/.style={
		line width=1.5pt, 
		square, 
		draw, 
		font=\sffamily,
		rounded corners, 
		inner sep=2pt
	},
	edge/.style={
		-stealth,
		shorten >=1pt,
		shorten <=1pt,
		line width=1.5pt
	},
	hyperedge/.style={
		Round Cap-{Triangle[length=3mm, width=2mm]},
		line width=5pt
	}
	]
	\node[main node] (v0) at (-1,0) {$0$};
	\node[main node] (v1) at (0,0) {$1$};
	\node[main node] (v2) at (1,0) {$2$};
	
	\node[second node] (w0) at (0,2) {$0$};
	\node[second node] (w1) at (0,-2) {$1$};
	
	\path[line width=1.5pt]
	(v0) edge [edge, in=210, out=240, looseness=8] (v0)
	(v1) edge [edge, in=210, out=240, looseness=8] (v1)
	(v2) edge [edge, in=300, out=330, looseness=8] (v2)

	(w0) edge [edge, uhhsteingrau, loop above, looseness=8] (w0)
	(w1) edge [edge, uhhsteingrau, loop below, looseness=8] (w1)
	;

        \path[line width=1.5pt]
	(v0) edge [edge, uhhdarkgrey, in=165, out=195, looseness=8] (v0)
	(v0) edge [edge, grey, in=120, out=150, looseness=8] (v0)
	
	(v1) edge [edge, uhhdarkgrey, in=30, out=60, looseness=8] (v1)
	(v0) edge [edge, grey] (v1)
	
	(v1) edge [edge, uhhdarkgrey] (v2)
	(v2) edge [edge, grey, in=30, out=60, looseness=8] (v2)
	;
	
	\coordinate (he1-1) at (-.75,1.5);
	\coordinate (he1-2) at (0,1.1);
	\coordinate (he1-3) at (.75,1.5);
	\coordinate (he1-4) at (-.75,-1.5);
	\coordinate (he1-6) at (0,-1.1);
	\coordinate (he1-5) at (.75,-1.5);
	
	\path[line width=1pt]
	(v0) edge [red, in=south west, out=north] (he1-1)
	(v1) edge [blue, in=south west, out=north west] (he1-1)
	
	(v1) edge [red, in=south, out=north] (he1-2)
	(v2) edge [blue, in=south, out=north west] (he1-2)
	
	(v2) edge [red, in=south east, out=north] (he1-3)
	(v0) edge [blue, in=south east, out=north east] (he1-3)
	
	(v0) edge [red, in=north west, out=south] (he1-4)
	(v2) edge [blue, in=north west, out=south west] (he1-4)
	
	(v2) edge [red, in=north east, out=south] (he1-5)
	(v1) edge [blue, in=north east, out=south east] (he1-5)
	
	(v1) edge [red, in=north, out=south] (he1-6)
	(v0) edge [blue, in=north, out=south east] (he1-6)
	;

	\draw[hyperedge, mix] (he1-1) -- (w0);
	\draw[hyperedge, mix] (he1-2) -- (w0);
	\draw[hyperedge, mix] (he1-3) -- (w0);
	\draw[hyperedge, mix] (he1-4) -- (w1);
	\draw[hyperedge, mix] (he1-5) -- (w1);
	\draw[hyperedge, mix] (he1-6) -- (w1);
	
\end{tikzpicture}%
		}%
		\caption{The hypernetwork that realizes the ODE system \eqref{vectorfieldex1}.}
		\label{fig:1st_ex}
	\end{subfigure}
	\hfill
	\begin{subfigure}{.6\textwidth}
		\includegraphics[width=\linewidth]{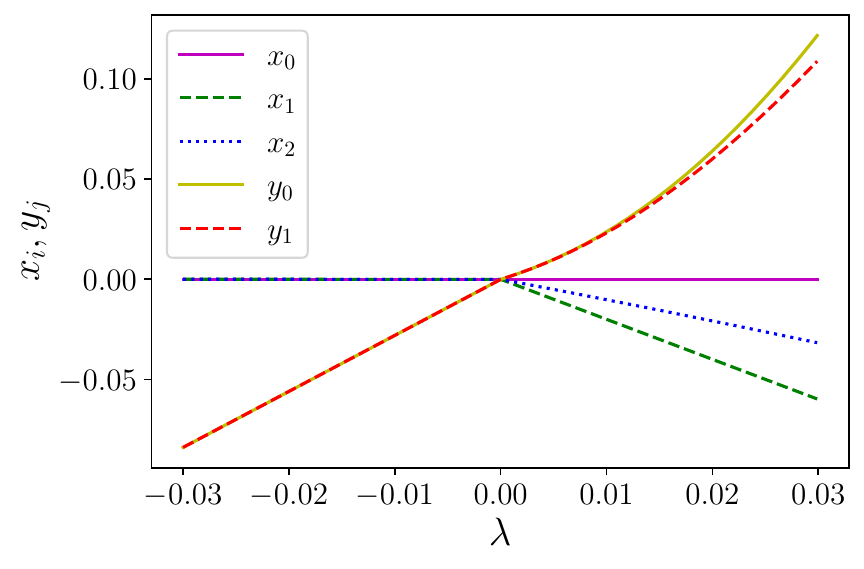}%
		\caption{Numerically obtained bifurcation diagram for equations of the form \eqref{vectorfieldex1}.
			Figure taken from \cite{vonderGracht.2023}.}
		\label{firstnumintro}
	\end{subfigure}
	\caption{A hypernetwork that supports an unusual ``reluctant''  synchrony  breaking local steady-state branch. The round nodes in \Cref{fig:1st_ex} correspond to the $x$-variables in \eqref{vectorfieldex1}, and the square nodes to the $y$-variables.
		The numerically obtained bifurcation  branch in \Cref{firstnumintro} satisfies $y_0=y_1$ for $\lambda<0$ and $y_0\neq y_1$ for $\lambda>0$. We will prove in this paper that $y_0-y_1 \sim \lambda^3$ for $\lambda>0$.}
\end{figure}

{\it Structure of the article.} In  \Cref{ex:leadingreluctant}, we illustrate reluctant synchrony breaking by studying the example presented in \cite{vonderGracht.2023} in more detail. In \Cref{Preliminaries}, we summarize the theoretical findings of \cite{vonderGracht.2023}, in which we defined a class of dynamical systems on hypernetworks, studied their robust synchrony and balanced partitions, and introduced the so-called augmented hypernetwork. In \Cref{sec:reluctant}, we prove our main result, \Cref{mainbifurcation}, which states  that reluctant synchrony breaking occurs generically in augmented hypernetworks, and  which provides a formula for the order of reluctancy of the synchrony breaking steady state branch. We also apply the theorem to the example discussed in this introduction. In \Cref{sec:moreexamples}, the main theorem is illustrated by three more examples. A discussion of our results is presented in \Cref{sec:discussion}.

\section{A first example}
\label{ex:leadingreluctant}
We now provide more details on the example that was briefly discussed in the introduction, and that was  introduced and studied numerically in \cite{vonderGracht.2023}.
As mentioned above,  this example concerns  the hypernetwork shown in \Cref{fig:1st_ex}.
In \cite{vonderGracht.2023}, we  introduced the class of {\it admissible ODEs} of a hypernetwork (see also \Cref{Preliminaries}). The admissible ODEs associated to the hypernetwork in \Cref{fig:1st_ex} are all the ODEs of the form
\begin{equation}\label{vectorfieldex1}
	\begin{array}{ll}
		\dot x_0 & = G(x_0, \textcolor{uhhdarkgrey}{x_0}, \textcolor{grey}{x_0})\, , \\
		\dot x_1 & = G(x_1, \textcolor{uhhdarkgrey}{x_1}, \textcolor{grey}{x_0}) \, , \\
		\dot x_2 & = G(x_2, \textcolor{uhhdarkgrey}{x_1}, \textcolor{grey}{x_2})\, , \\
		\dot y_0 & = F(\textcolor{uhhsteingrau}{y_0}, \textcolor{mix}{(} \textcolor{red}{x_0}, \textcolor{blue}{x_1} \textcolor{mix}{)}, \textcolor{mix}{(} \textcolor{red}{x_1}, \textcolor{blue}{x_2} \textcolor{mix}{)}, \textcolor{mix}{(} \textcolor{red}{x_2}, \textcolor{blue}{x_0} \textcolor{mix}{)} )\, ,\\
		\dot y_1 & = F(\textcolor{uhhsteingrau}{y_1}, \textcolor{mix}{(} \textcolor{red}{x_0}, \textcolor{blue}{x_2} \textcolor{mix}{)}, \textcolor{mix}{(} \textcolor{red}{x_1}, \textcolor{blue}{x_0} \textcolor{mix}{)}, \textcolor{mix}{(} \textcolor{red}{x_2}, \textcolor{blue}{x_1} \textcolor{mix}{)} )\, ,\\
	\end{array}
\end{equation}
for certain smooth functions $F$ and $G$. We assume for now that the variables $x_i, y_j$ take values in $\R$, so that $F: \R \times \R^2\times \R^2\times \R^2 \to \R$ and $G:\R\times \R\times \R \to \R$. 
The brackets in $F$ serve to distinguish the two-dimensional inputs from hyperedges of order two (the purple arrows in the \Cref{fig:1st_ex}).   
The assumption that these hyperedges are   identical, translates into the requirement that $F$ is invariant under all permutations of these  pairs of variables. 
That is, we require that for all $Y, X_1, \dots, X_6 \in \R$, 
\begin{equation}
	\label{exampleFhype}
	\begin{split}
		& F(Y, (X_0, X_1), (X_2, X_3), (X_4, X_5)) =  \\
		&  F(Y, (X_2, X_3), (X_0, X_1), (X_4, X_5)  ) = \\ 
		&  F(Y,(X_0, X_1), (X_4, X_5), (X_2, X_3) )\, .
	\end{split}
\end{equation}

\noindent To study bifurcations within this class of admissible ODEs, we parameterize the response functions $F$ and $G$ in  \eqref{vectorfieldex1} by a scalar variable $\lambda$ taking values in some open neighborhood $\Omega \subseteq \R$ of the origin. 
This gives the one-parameter family of hypernetwork admissible ODEs
\begin{equation}\label{systembifwithpara}
	\begin{array}{ll}
		\dot x_0 & = G(x_0, \textcolor{uhhdarkgrey}{x_0}, \textcolor{grey}{x_0}; \lambda)\, , \\
		\dot x_1 & = G(x_1, \textcolor{uhhdarkgrey}{x_1}, \textcolor{grey}{x_0}; \lambda) \, , \\
		\dot x_2 & = G(x_2, \textcolor{uhhdarkgrey}{x_1}, \textcolor{grey}{x_2}; \lambda)\, , \\
		\dot y_0 & = F(y_0, \textcolor{mix}{(} \textcolor{red}{x_0}, \textcolor{blue}{x_1} \textcolor{mix}{)}, \textcolor{mix}{(} \textcolor{red}{x_1}, \textcolor{blue}{x_2} \textcolor{mix}{)}, \textcolor{mix}{(} \textcolor{red}{x_2}, \textcolor{blue}{x_0} \textcolor{mix}{)}; \lambda)\, ,\\
		\dot y_1 & = F(y_1, \textcolor{mix}{(} \textcolor{red}{x_0}, \textcolor{blue}{x_2} \textcolor{mix}{)}, \textcolor{mix}{(} \textcolor{red}{x_1}, \textcolor{blue}{x_0} \textcolor{mix}{)}, \textcolor{mix}{(} \textcolor{red}{x_2}, \textcolor{blue}{x_1} \textcolor{mix}{)}; \lambda)\, .\\
	\end{array}
\end{equation}
To guarantee that the system \eqref{systembifwithpara} is admissible for every fixed value of $\lambda \in \Omega$, we  assume $F$ satisfies Equation \eqref{exampleFhype} for any fixed value of $\lambda$.
Let us in addition assume 
\begin{equation*}
	F(0,(0,0), (0,0), (0,0);0) =  G(0,0,0;0) = 0\, ,
\end{equation*}
so that the system \eqref{systembifwithpara} has a steady-state point at the origin for $\lambda = 0$. We may study the persistence of this steady-state by investigating the Jacobian at the origin of the system for $\lambda = 0$. 
To this end, we write
\begin{equation}\label{expressionsforFnnfG}
	\begin{split}
		& F(Y, (X_0, X_1), (X_2, X_3), (X_4, X_5); \lambda )  \\   
		& = aY + bX_0 + cX_1 + bX_2 + cX_3 + bX_4 + cX_5 + d\lambda  + \mathcal{O}(\|(Y,X_0, \dots, X_5;\lambda)\|^2) \, ,
	\end{split}
\end{equation}
and
\begin{equation}\label{expressionsforFnnfG2}
	G(X_0,X_1,X_2;\lambda) = AX_0+BX_1+CX_2 + \mathcal{O}(|\lambda| + \|(X_0,X_1,X_2)\|^2)\, ,
\end{equation}
with $a, \dots, d, A, B, C \in \R$, to specify the linear terms. The multiple occurrence of the terms $b$ and $c$ in \eqref{expressionsforFnnfG} is due to the fact that according to \eqref{exampleFhype},  $F$ depends in the same way on the tuple $(X_0,X_1)$ as it does on $(X_2,X_3)$ and $(X_4,X_5)$.
In terms of these coefficients, the Jacobian matrix of the right hand side of \Cref{systembifwithpara} at $(x; \lambda)=(0,0)$ with respect to the spatial variables $x=(x_0,x_1,x_2,y_0,y_1)$ is
$$\left( \begin{array}{ccccc} 
	A+\textcolor{uhhdarkgrey}{B}  + \textcolor{grey}{C}& 0 & 0 & 0 & 0 \\
	\textcolor{grey}{C} & A+\textcolor{uhhdarkgrey}{B} & 0 & 0 & 0 \\
	0 & \textcolor{uhhdarkgrey}{B} & A+\textcolor{grey}{C} & 0 & 0 \\ 
	\textcolor{red}{b}+\textcolor{blue}{c} & \textcolor{red}{b}+\textcolor{blue}{c}  & \textcolor{red}{b}+\textcolor{blue}{c} & a & 0 \\
	\textcolor{red}{b}+\textcolor{blue}{c} & \textcolor{red}{b}+\textcolor{blue}{c} & \textcolor{red}{b}+\textcolor{blue}{c} & 0  & a 
\end{array} 
\right)\, .
$$
The eigenvalues of this Jacobian are $A+B+C$, $A+B$ and $A+C$ (all with multiplicity $1$), and $a$ (with algebraic and geometric multiplicity $2$). 
To allow for a steady state bifurcation to occur at $\lambda=0$, we consider the case $A+B = 0$. We moreover assume the generic conditions $a, A+B+C, A+C \not= 0$ to hold.

We claim that as $\lambda$ is varied near $0$, two branches of steady states will generically emerge from the origin. These can be found by first focusing on the subnetwork given by the three nodes of the same type. That is, we first solve   
\begin{equation*}
	\begin{split}
		&G(x_0, {x_0}, {x_0}; \lambda) = 0 \, ,\\ \nonumber
		&G(x_1, {x_1}, {x_0}; \lambda)= 0\, ,\\ \nonumber
		&G(x_2, {x_1}, {x_2}; \lambda) = 0\, .
	\end{split}
\end{equation*}
A direct calculation shows that, for generic values of the first and second degree Taylor coefficients of $G$, one of the steady state branches is locally given by 
\begin{equation}\label{firstbranchs}
	x_0(\lambda) = x_1(\lambda) = x_2(\lambda) = x(\lambda) = D_0\lambda + \mathcal{O}(|\lambda|^2) \, ,
\end{equation}
while another branch is given by 
\begin{equation}\label{secondbranchs}
	x_0(\lambda) = D_0\lambda + \mathcal{O}(|\lambda|^2) \, , \ x_1(\lambda) = D_1\lambda + \mathcal{O}(|\lambda|^2) \, , 
	x_2(\lambda) = D_2\lambda + \mathcal{O}(|\lambda|^2)\, ,
\end{equation}
for certain non-zero and mutually  distinct $D_0, D_1, D_2 \in \R$. 
For our choices of parameters, no further branches exist.
We omit the computation of these branches. 
For a  detailed exposition on how to compute steady state bifurcation branches in  so-called \emph{feedforward networks} (i.e. networks with no loops other than self-loops), we refer to \cite{vondergracht2022}. 

We now turn to computing the values of $y_0$ and $y_1$ along the bifurcation branches. 
We start by looking at the first branch, given by \Cref{firstbranchs}. Restricted to this branch, the steady state equation $\dot{y}_0 = 0$ becomes 
\begin{equation}\label{y0is0}
	F(y_0, (x(\lambda), x(\lambda)), (x(\lambda), x(\lambda)), (x(\lambda), x(\lambda)); \lambda ) = 0\, .
\end{equation}
Combining   \eqref{expressionsforFnnfG} and \eqref{firstbranchs} this can be expanded as
\begin{equation*}
	ay_0 + (3D_0(b+c) + d)\lambda 
	+ \mathcal{O}(\|(y_0;\lambda)\|^2) = 0\, ,
\end{equation*}
which by the implicit function theorem has a unique solution given by
\begin{equation*}
	y_0(\lambda) = \frac{-3D_0(b+c) - d}{a}\lambda  + \mathcal{O}(|\lambda|^2)\, .
\end{equation*}
Setting $\dot{y}_1 = 0$ gives precisely the same  equation to solve as \eqref{y0is0}, but with $y_0$ replaced by $y_1$. 
Hence, we find $y_0(\lambda) = y_1(\lambda)$ along this first branch, which we will therefore refer to as the \emph{synchronous branch} of system \eqref{systembifwithpara}.

We now turn to the second branch of steady states, of which the asymptotics of the $x$-variables is given by \Cref{secondbranchs}.
Combining \eqref{expressionsforFnnfG} with \eqref{secondbranchs}, we find that $\dot{y}_0 = 0$ is equivalent to
\begin{equation*}
	ay_0 + ((b+c)(D_0 + D_1 + D_2) + d)\lambda 
	+ \mathcal{O}(\|(y_0;\lambda)\|^2) = 0\, .
\end{equation*}
It follows again from the implicit function theorem that locally precisely one solution exists, given by 
\begin{equation}\label{y0non-synns}
	y_0(\lambda) = \frac{-(b+c)(D_0 + D_1 + D_2) - d}{a}\lambda  + \mathcal{O}(|\lambda|^2)\, .
\end{equation}
In exactly the same way, we find that $\dot{y}_1 = 0$ is solved for by 
\begin{equation}\label{y1non-synns}
	y_1(\lambda) = \frac{-(b+c)(D_0 + D_1 + D_2) - d}{a}\lambda  + \mathcal{O}(|\lambda|^2)\, .
\end{equation}
Note that these expressions for $y_0(\lambda)$ and $y_1(\lambda)$ agree up to first order in $\lambda$.
However, unlike for the synchronous branch, there is no reason to conclude that $y_0(\lambda) = y_1(\lambda)$ along this branch, as the equations for $\dot{y}_0$ and $\dot{y}_1$ in \eqref{systembifwithpara} are  different for distinct $x_0, x_1$ and $x_2$:
\begin{align*}
	F(y_0, (x_0(\lambda), x_1(\lambda)), (x_1(\lambda), x_2(\lambda)), (x_2(\lambda), x_0(\lambda)); \lambda ) &= 0\, ,\\
	F(y_1, (x_0(\lambda), x_2(\lambda)), (x_1(\lambda), x_0(\lambda)), (x_2(\lambda), x_1(\lambda)); \lambda ) &= 0\, .
\end{align*}
We will refer to the branch of steady states given by \Cref{secondbranchs,y0non-synns,y1non-synns} as the \emph{reluctant  branch} of \eqref{systembifwithpara}. 

\begin{figure}
	\centering
	\begin{subfigure}{0.47\linewidth}
		\centering
		\includegraphics[width=\linewidth]{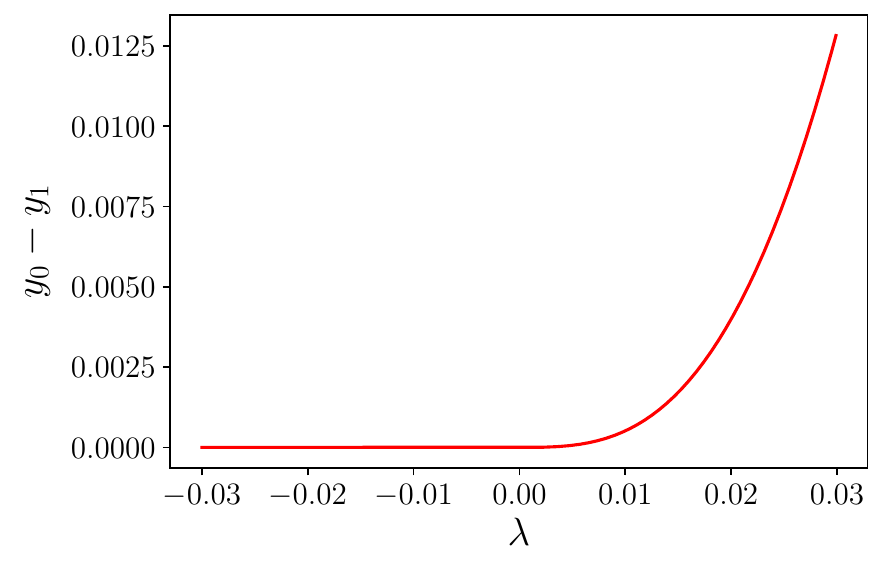}
		\caption{The difference between the $y$-nodes along the steady state branches.}
		\label{firstnum-b}
	\end{subfigure}\hfill
	\begin{subfigure}{0.47\linewidth}
		\centering
		\includegraphics[width=\linewidth]{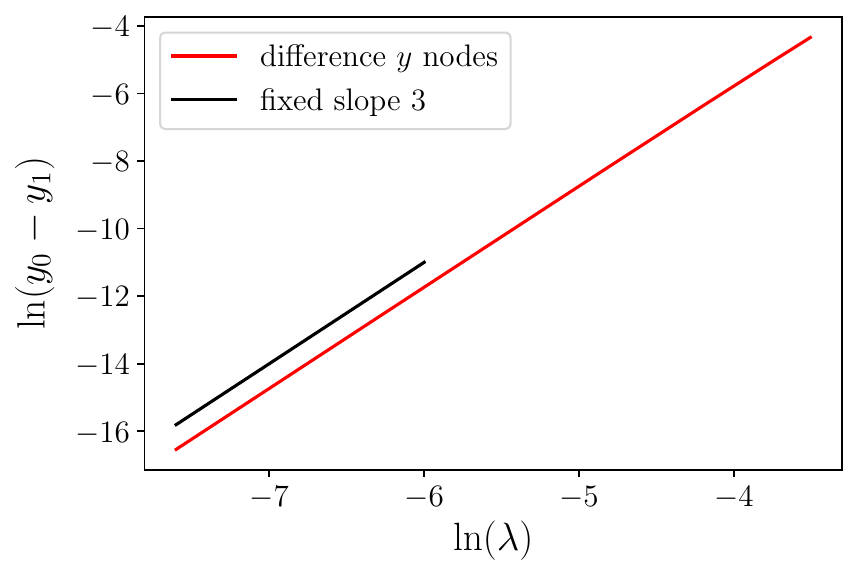}
		\caption{A log–log plot of the difference between the $y$-nodes, for $\lambda > 0$.}
		\label{firstnum-c}
	\end{subfigure}
	\caption{More details for the steady state branches depicted in \Cref{firstnumintro}. The black line segment in the log-log plot has slope $3$ and was added to show that $y_0(\lambda) - y_1(\lambda) \sim \lambda^3$. 
		Figures taken from \cite{vonderGracht.2023}. }
	\label{firstnum}
\end{figure}

\Cref{firstnumintro} demonstrates numerically that the reluctant branch is truly non-synchronous. The figure was taken from \cite{vonderGracht.2023}, and it shows a numerically obtained plot of the asymptotically stable bifurcation branches that emerge in a bifurcation in  a particular realisation of system \eqref{systembifwithpara}, namely for the choices 
\begin{align}
	& \nonumber 
	G(X_0,X_1,X_2; \lambda) =  -X_0+X_1-X_2+8\lambda X_0+4X_0^2 \ \text{ and } \\ 
	& \nonumber 
	F(Y, (X_0, X_1), (X_2, X_3), (X_4, X_5); \lambda ) =  -5Y + 14 \lambda 
	\\ \nonumber
	&-h(10X_0-12X_1)-h(10X_2-12X_3)-h(10X_4-12X_5)   \,  
\end{align}
in which
\begin{equation}
	h(x) = \sin(x)+\cos(x)-1  \, .
\end{equation}
It is known from \cite{vondergracht2022} that for this choice of $G$, the steady state branch $x(\lambda)$ is stable inside the subnetwork given by the three nodes of the same type. The negative coefficient in front of the linear $Y$-term in $F$ implies that it is stable in the $y$-directions as well --- so that it can easily be found  numerically. \Cref{firstnumintro} shows the synchronous branch for $\lambda<0$ and the reluctant branch for $\lambda>0$ -- 
indeed, $y_0$ and $y_1$ agree for $\lambda<0$, 
and quite clearly do not for $\lambda>0$. This becomes even more  visible in \Cref{firstnum-b}, which shows $y_0-y_1$ as a function of $\lambda$. The logarithmic plot in \Cref{firstnum-c} suggests that $y_0(\lambda) - y_1(\lambda) \sim \lambda^3$.
In \Cref{sec:reluctant}, we prove that this is truly the case.

In this article, we rigorously prove the existence of reluctant steady state branches in bifurcations in a large class of hypernetwork systems. Specifically, we will show that  reluctant bifurcation branches appear in generic one-parameter local synchrony breaking steady state bifurcations in the admissible ODEs for such hypernetworks. 
We also provide a formula for the order in $\lambda$ with which the ``reluctant nodes'' separate.

\section{Preliminaries}
\label{Preliminaries}
In this section, we briefly introduce hypernetwork dynamical systems and summarize results obtained in \cite{vonderGracht.2023}.
In fact, in \cite{vonderGracht.2023} we define a  hypernetwork to be a collection ${\rm \bf N} = (V, H, s, t)$ consisting of a finite set of \emph{vertices} or \emph{nodes} $V$, a finite set of \emph{hyperedges} $H$, and source and target maps $s$ and $t$ defined on $H$. 
Given an edge $h \in H$, its target $t(h) \in V$ is a single vertex, whereas the source $s(h) = (s_1(h), \dots, s_{k_h}(h)) \in V^{k_h}$ is an ordered ${k_h}$-tuple of vertices. 
The number ${k_h}>0$ depends on the hyperedge $h$, and is called its \emph{order}. 
To avoid cluttered notation though, we often suppress the dependence of ${k_h}$ on $h$ when it is clear from context, and simply write $k$.
Note that the $k_h$ vertices in $s(h)$ are not required to be distinct.
The order of a hypernetwork ${\rm \bf N}$ is then defined as the maximum of the orders of its hyperedges, so that the hypernetworks of order $1$ are precisely the classical (dyadic) networks.

In addition to the data that is explicitly given in ${\rm \bf N} = (V, H, s, t)$, we also specify equivalence relations on both the nodes $V$ and the hyperedges $H$. We typically refer to both as the \emph{color} or \emph{type} relation. 
The reason that these relations are not specified in ${\rm \bf N}$ is because they will apply to all hypernetworks at once.
That is, it will make sense for two nodes in different hypernetworks to have the same color, and likewise for multiple hyperedges across different hypernetworks.
This allows us to define node- and hyperedge-type preserving maps between different hypernetworks, which in turn give rise to semi-conjugacies between the dynamics, see \cite{vonderGracht.2023} for more details.
Intuitively, this color-relation conveys whether two nodes correspond to comparable or incomparable agents in a real-world system modelled by the  
hypernetwork, and similarly whether or not two hyperedges specify the same influence.
As is suggested by this interpretation, the vertex- and hyperedge-types have to satisfy certain consistency conditions. These are:
\begin{enumerate}
	\item if two nodes $v_0$ and $v_1$ are of the same type, then there exists a hyperedge-type preserving bijection between the set of hyperedges targeting $v_0$ and those targeting $v_1$;
	\item two hyperedges $h_0, h_1$ of the same type have the same order $k$, and for each $i \in \{1, \dots, k\}$ the nodes $s_i(h_0)$ and $s_i(h_1)$ are of the same vertex-type.
\end{enumerate}

\noindent These conditions allow us to define dynamical systems with the interaction structure of the given hypernetwork. Such dynamical systems are specified by so-called \emph{admissible} vector fields. To introduce these, we fix an internal phase space $\R^{n_v}$ for each node $v \in V$, which will be identical for nodes of the same type.
Each node $v$ is given a state variable $x_v \in \R^{n_v}$. The  \emph{total phase space} of the hypernetwork dynamical system describes the states of all these variables, and  is thus the direct sum $ \bigoplus_{v\in V} \mathbb{R}^{n_v}$. 
We also specify, for each hyperedge $h \in H$, the vector of its source variables \[ {\bf x}_{s(h)} = \left(x_{s_1(h)}, \ldots, x_{s_{k_h}(h)}\right) \, \in \, \bigoplus_{i = 1} ^ {k_h} \mathbb{R}^{n_{s_i(h)}}. \]
Next, we choose for each node $v \in V$ its \emph{response function} 
\begin{align}
	F_v \colon   \bigoplus_{h\, :\, t(h)=v} \bigoplus_{i = 1} ^ {k_h} \, \, \mathbb{R}^{n_{s_i(h)}} \to \R^{n_v}\, .
\end{align} 
These functions must satisfy certain conditions that reflect our intuitive idea that hyperedges of the same type encode identical influence, as well as the notion that nodes of the same type respond to their input in the same way.
In words, we require that the variables of identical-type hyperedges may be freely interchanged in $F_v$, as well as that $F_v$ and $F_w$ are the same when the nodes $v$ and $w$ are of the same type, after an appropriate identification of their domains.
We may capture both requirements in one succinct condition as follows: given nodes $v$ and $w$ of the same type, for any hyperedge-type preserving bijection $\alpha: t^{-1}(v) \to t^{-1}(w)$ we have 
\begin{equation}
	\label{eq:bijection_symmetry}
	F_{w}\left( \bigoplus_{t(h_2)=w} {\bf x}_{s(h_2)} \right) = F_{v}\left(\bigoplus_{t(h_1)=v} {\bf x}_{s(\alpha(h_1))} \right) \, ,
\end{equation}
for all $x = \bigoplus_{v \in V} x_v \in \bigoplus_{v\in V} \mathbb{R}^{n_v}$. 
Recall that at least one such $\alpha$ exists when $v$ and $w$ are of the same type.
Finally, we define the {\it hypernetwork admissible vector field} $$f^{\bf N}: \bigoplus_{v\in V} \mathbb{R}^{n_v} \to  \bigoplus_{v\in V} \mathbb{R}^{n_v}\, ,$$
on the total phase space, to be given component-wise by 
$$f^{\bf N}_v(x) = F_v\left(  \bigoplus_{h\, :\, t(h)=v} {\bf x}_{s(h)}  \right)\,  $$
for all $v \in V$ and $x \in \bigoplus_{v\in V} \mathbb{R}^{n_v}$.

To study synchronisation in hypernetwork dynamical systems, we define a {\it polysynchrony subspace} to be a subspace of the total phase space of a hypernetwork that is determined by equality of clusters of node variables. More precisely, for any partition  $P=\{V_1,\dotsc,V_C\}$ of the nodes $V$ of a hypernetwork, one can define the  polysynchrony subspace
\[ \operatorname{Syn}_P=\{ x_{v}=x_{w} \text{ when } v \text{ and } w \text{ are in the same element of } P\}. \]
Such a polysynchrony subspace is called {\it robust}  if it is invariant under the flow  of any admissible vector field, that is, when $f^{\bf N}(\operatorname{Syn}_P)\subseteq \operatorname{Syn}_P$ for every admissible vector field $f^{\bf N}$.   
It was shown in \cite{vonderGracht.2023} that $\operatorname{Syn}_P$ is robust if and only if $P$ is \emph{balanced}, meaning that the partition is ``consistent with the hypernetwork structure''. For the precise definition of a balanced partition, we refer to \cite{vonderGracht.2023}. It was also shown in \cite{vonderGracht.2023} that a partition $P$ is balanced (and hence $\operatorname{Syn}_P$ is robust) if and only if $\operatorname{Syn}_P$ is  invariant under all polynomial admissible vector fields of degree at most $\frac{k(k+1)}{2}$, where $k$ is the order of the hypernetwork. 
\begin{ex}
	One can show that the hypernetwork discussed in \Cref{ex:leadingreluctant} and depicted in \Cref{fig:1st_ex} has four robust synchrony subspaces (apart from the total phase space itself), namely $\{x_0=x_1\}, \{x_0=x_1=x_2\}$, $\{x_0=x_1 \text{ and } y_0=y_1\}$,  and $\{x_0=x_1=x_2 \text{ and } y_0=y_1\}$.
\end{ex}

\noindent 
Of particular interest in this paper are so-called \emph{augmented hypernetworks}, also introduced in \cite{vonderGracht.2023}. 
Their definition involves the symmetric group on $k+1$ elements, denoted by $S_{k+1}$, which acts on the ordered set $\{0, \dots, k\}$ by permutations.
We  denote by $S_{k+1}^0$ and $S_{k+1}^1$ the subsets of even and odd permutations, respectively, and denote by $\sgn(\sigma) \in \{0,1\}$ the sign of a permutation ${\sigma \in S_{k+1}^{\sgn(\sigma)}}$.

\begin{defi}[Definition 5.1 in \cite{vonderGracht.2023}]\label{def:augmentedhere}
	Let ${\rm \bf N}$ be a hypernetwork with $k+1 \geq 3$ nodes $v_0, \dots, v_k$ of the same type. We define the \emph{augmented hypernetwork with core ${\rm \bf N}$}, denoted by $\NN^\diamondsuit$, as the hypernetwork obtained by adding two additional nodes $w_0, w_1$, one self-loop for each $w_i$ and $(k+1)!$ new hyperedges to  ${\rm \bf N}$. 
	The new nodes are of the same type, which differs from that of the $v_i$. 
	Likewise, we construct a new hyperedge-type which we assign to the $(k+1)!$ additional hyperedges.
	Necessarily the self-loops on the new nodes are of a same, new type too.
	The $(k+1)!$ new hyperedges are indexed by the symmetric group $S_{k+1}$, so that we may denote them by $h_{\sigma}$ for $\sigma \in S_{k+1}$. 
	We define their sources and targets by 
	\begin{equation}
		t(h_{\sigma}) = w_{\sgn(\sigma)} \quad \text{and} \quad
		s(h_{\sigma}) = (v_{\sigma(1)}, \dots, v_{\sigma(k)}) \, ,
	\end{equation}
	where $S_{k+1}$ acts on the ordered set $\{0, \dots, k\}$.
	Note that $v_{\sigma(0)} \in\{v_0, \dots, v_k\}$ is therefore the only $v$-node not in the source of $h_{\sigma}$, and that these hyperedges all have order $k$. 
\end{defi}

\begin{ex}
	The hypernetwork discussed in \Cref{ex:leadingreluctant} and depicted in \Cref{fig:1st_ex}, is an example of an augmented hypernetwork. Here the core consists of the $k+1 = 3$ circular nodes in the center and the arrows between them.  The core in fact forms a classical (dyadic) network. The two ``added'' nodes are depicted as the square ones. 
	\Cref{vectorfieldex1} gives the form of a general admissible vector field for this augmented hypernetwork. 
	Recall that the response function $F$ in \Cref{vectorfieldex1} is invariant under permutations of the three pairs of inputs, which reflects that the six ``added'' hyperedges are all of the same type. 
\end{ex}

\noindent By assumption, all nodes in the  core $\NN$ of an augmented hypernetwork are of the same type, so that $\NNd$ has precisely two node-types.
This means that two response functions are required to describe an admissible vector field $f^{\NNd}$  for $\NN^\diamondsuit$. 
We will usually denote these by $F$ and $G$, where $G$ is used for the nodes in $\NN$ and $F$ for the two additional nodes. 
Likewise, we see that the total phase space is determined by two vector spaces: one for the internal dynamics of the $v$-nodes, $\R^{n_v}$, and one for that of the $w$-nodes, $\R^{n_w}$.
We will later set both equal to $\R$.
Note that $F$ takes one argument from $\R^{n_w}$, corresponding to the self-loop, and $(k+1)!/2$ entries from $\bigoplus^k \R^{n_v}$ for the remaining hyperedges.
As these latter hyperedges are indexed by (half of) the symmetric group, we may see the response function as
$$  F \colon\R^{n_w} \oplus \hspace{-3pt}\bigoplus_{\sigma \in S^0_{k+1}}  \hspace{-5pt} \bigg(\bigoplus^k \R^{n_v}\bigg) \, \to\, \R^{n_w}\, ,$$
with the property that the $(k+1)!/2$ entries with values in $\bigoplus^k \R^{n_v}$ may be freely interchanged.
Note that we simply index these entries by $S^0_{k+1}$ to emphasize that they correspond to the hyperedges that are indexed by (part of) the symmetric group. 
We could have also used $S^1_{k+1}$ and, because we may freely interchange these entries, we do not have to give an explicit identification between $S^0_{k+1}$ and $S^1_{k+1}$.
In particular, in any augmented hypernetwork the dynamics of the $w$-nodes may simply be written as
\begin{equation}
	\dot{y}_{0} = F\left({y}_{0}, \bigoplus_{\sigma \in S^0_{k+1}} {\bf x}_{\sigma} \right) \text{ and }\, \dot{y}_{1} = F\left({y}_{1}, \bigoplus_{\sigma \in S^1_{k+1}} {\bf x}_{\sigma} \right)\, ,
\end{equation}
where ${\bf x}_{\sigma} = (x_{\sigma(1)}, \dots, x_{\sigma(k)})$, and where we write $y_i$ for the state of node $w_i$ and $x_j$ for the state of node $v_j$.

\begin{remk}\label{moregeneralaug}
	We may generalise the definition of an augmented hypernetwork by connecting the auxiliary nodes $w_0$ and $w_1$ only to a subset $\mathcal{S}$ of nodes of the same type of the core, consisting of $k+1$ nodes of the same type. In this case, the nodes in the core that are not elements of $\mathcal{S}$ may in fact be of different type.
	We can then add $(k+1)!$ hyperedges of order $k$, precisely as in \Cref{def:augmentedhere}, but now with all source nodes in $\mathcal{S}$.
	It will be clear that our results also hold for such hypernetworks, but to avoid a cluttered exposition we will mostly work with \Cref{def:augmentedhere}. 
	See also \Cref{remk:firstclosing}. 
\end{remk}

\section{Reluctant synchrony-breaking}
\label{sec:reluctant}
We now show that the reluctant synchrony breaking observed in \Cref{ex:leadingreluctant} is not a peculiarity of systems of the form \eqref{systembifwithpara}, but can occur in any augmented hypernetwork. 
In fact, we shall give natural conditions on the core $\NN$ that guarantee that reluctant synchrony  breaking occurs generically in the augmented hypernetwork $\NN^{\diamondsuit}$. Moreover, in \Cref{mainbifurcation} below we  give a precise expression for the degree (in the bifurcation parameter) at which the reluctant synchrony breaking occurs. 
We start by introducing some useful notation and conventions.

As it is sometimes convenient to make explicit the dependence of an admissible vector field on its response functions, we will often write $f^{\NNd}_{(F,G)}$ and $f^{\NN}_{(G)}$ for the admissible vector fields of $\NNd$ and $\NN$, respectively. 
Furthermore, because in this section we are mainly interested in bifurcations, we will often use $f^{\NN}$ (and  $f^{\NNd}$, $f^{\NNd}_{(F,G)}$ etc.) to denote parameter families of admissible vector fields.
This means $f^{\NN}$ is an admissible vector field for any fixed value of the bifurcation parameter, as in \Cref{ex:leadingreluctant}.

Throughout this section we will investigate asymptotics and power series in $\lambda$ for bifurcation branches. Some of these might involve fractional powers of $\lambda$, meaning that such branches are only defined for positive or negative values of $\lambda$. To keep this section as readable as possible, we assume from here on out that all branches are defined for positive values of $\lambda$, so that we may always write $\lambda^p$ for any  power $p \geq 0$. 
The corresponding results for negative values of $\lambda$ follow easily by redefining $\lambda$ as $-\lambda$. For a function $X=X(\lambda)$ defined for small positive values of $\lambda$,  we will frequently use the abbreviation $X(\lambda)\sim\lambda^p$ to denote that  there exists a constant $A\ne 0$  for which $X(\lambda)=A\lambda^p +  \text{ ``higher order terms in } \lambda$''.

\begin{defi}\label{def:orderandorderofasyn}
	Let $\NN$ be a hypernetwork with $n$ nodes and denote by $f^{\NN} \colon \R^n \times \Omega \rightarrow \R^n$ a one-parameter family of admissible vector fields for $\NN$, where each node has a one-dimensional internal phase space.
	In this paper, a locally defined branch of steady states $x(\lambda) = (x_1(\lambda), \dots, x_n(\lambda))$ for $f^{\NN}$ is called \emph{fully synchrony-breaking} if for all $i,j \in \{1, \dots, n\}$ with $i>j$ there exist numbers $p_{i,j}>0$
	such that
	\begin{equation}\label{differ4nces}
		x_i(\lambda) - x_j(\lambda) \sim \lambda^{p_{i,j}}\, .
	\end{equation}
	Given a fully synchrony-breaking steady state branch, we define its \emph{order of asynchrony} as the number
	\begin{equation}
		\bar{p} = \sum_{\substack{i,j = 1 \\ i>j}}^n p_{i,j}\, .
	\end{equation}
\end{defi}
\noindent In what follows, we will need to make some assumptions about the asymptotics and regularity of the branches. 
To avoid a detailed exposition, we instead use a formal ansatz. 
From here on out, we will always assume a fully synchrony-breaking branch $x(\lambda) = (x_1(\lambda), \dots, x_n(\lambda))$ to come with a finite set of numbers $\Upsilon \subseteq (0, \bar{p}\,]$ such that we may write
\begin{equation}\label{eq:formalpowerseries}
	x_i(\lambda) = \sum_{p \in \Upsilon} A_{i,p} \lambda^{p} + \mathcal{O}(|\lambda|^{\bar{p}+\epsilon})\, 
\end{equation}
for some $\epsilon > 0$, and with $A_{i,p} \in \R$. 
Note that these $A_{i,p}$ may very well vanish.
As we may write
\begin{align}\label{eq:formalpowerseriesdifff}
	x_i(\lambda)-x_j(\lambda) &= \sum_{p \in \Upsilon} (A_{i,p}-A_{j,p}) \lambda^{p} + \mathcal{O}(|\lambda|^{\bar{p}+\epsilon}) \\ \nonumber
	&= D_{i,j}\lambda^{p_{i,j}} + \text{ ``higher order terms''} \, ,
\end{align}
and because $p_{i,j} \leq \bar{p}$, we see that necessarily $p_{i,j} \in \Upsilon$ for all $i>j$. Note that the $D_{i,j}$ appearing in \eqref{eq:formalpowerseriesdifff} are all nonzero by assumption \eqref{differ4nces}.

Equation \eqref{eq:formalpowerseriesdifff} also shows that $A_{i,p_{i,j}}$ and $A_{j,p_{i,j}}$ cannot both be $0$, as this would imply $D_{i,j}=0$. 
As a result, we see that at least one of the components $x_i(\lambda)$ and $x_j(\lambda)$ grows as $\lambda^s$ with $s \leq p_{i,j}$. 
Given $k \in \{1, \dots, n\}$ let us denote by $s_k$ the order of $x_k(\lambda)$, so that $x_k(\lambda) \sim \lambda^{s_k}$. 
If $x_k(\lambda) = 0$ then we may set $s_k = \infty$.
What we have argued is simply that $\min(s_i, s_j) \leq p_{i,j}$.

A related quantity appearing in Theorem \ref{mainbifurcation} below will be
\begin{equation}
	\hat{p} = \min(s_1, \dots, s_n, 1) \, .
\end{equation}
It follows that
\begin{equation}\label{phatpbar}
	\hat{p} \leq \min(s_1, s_2) \leq p_{1,2} < \bar{p} \, ,
\end{equation}
as $p_{i,j} > 0$ for all $i,j$.

In the theorem below, we assume all admissible vector fields correspond to one-dimensional internal dynamics for each node. 
\begin{thr}\label{mainbifurcation}
	Let $\NNd$ be an augmented hypernetwork with core $\NN$, the latter consisting of $k+1 \ge 3$ nodes, and let $f^{\NN}_{(G)}: \R^{k+1} \times \Omega \rightarrow \R^{k+1}$ be a one-parameter family of admissible vector fields for $\NN$, corresponding to some response function $G$. 
	Assume $f^{\NN}_{(G)}$ admits a fully synchrony-breaking branch of steady states $x(\lambda) = (x_0(\lambda),   \dots, x_{k}(\lambda))$ with order of asynchrony $\bar{p}$.
	
	Then for a generic $\lambda$-dependent response function $F=F(\cdot ;\lambda)$, the system $f^{\NNd}_{(F,G)}: \R^{k+3} \times \Omega \rightarrow \R^{k+3}$ admits a steady state branch 
	\[ z(\lambda) = (x_0(\lambda), \dots, x_k(\lambda), y_0(\lambda), y_1(\lambda)) \]
	for which 
	\begin{equation*}
		y_0(\lambda), y_1(\lambda) \sim \lambda^{\hat{p}}\ \mbox{while} \ 
		y_0(\lambda) - y_1(\lambda) \sim  \lambda^{\bar{p}} \,  .
	\end{equation*} 
	The branch $z(\lambda)$ consists of \emph{asymptotically stable} equilibria for $f^{\NNd}_{(F,G)}$ if and only if $x(\lambda)$ consists of asymptotically stable equilibria for $f^{\NN}_{(G)}$ and $\frac{\partial F}{\partial y}(0;0)<0$.
\end{thr}

\noindent Note that, since $\hat{p} < \bar{p}$ by Equation \eqref{phatpbar}, we see that the branch found in Theorem \ref{mainbifurcation} is truly tangent to the space $\{y_0 = y_1\}$. In other words, the difference between the $y$-components grows significantly slower than the two $y$-values themselves.

\begin{remk}
	As $(f^{\NNd}_{(-F,G)})_{w} = -(f^{\NNd}_{(F,G)})_{w}$ for the two nodes $w$ outside the core, we see that $z(\lambda)$ is a branch of steady-states for $f^{\NNd}_{(F,G)}$, if and only if it is for $f^{\NNd}_{(-F,G)}$.
	The stability condition in \Cref{mainbifurcation} holds for either $f^{\NNd}_{(F,G)}$ or $f^{\NNd}_{(-F,G)}$. 
	Hence, each existing branch that is stable for the core is stable for either $f^{\NNd}_{(F,G)}$ or for $f^{\NNd}_{(-F,G)}$. This can be particularly useful if one wants to realise all branches (not just the stable ones) numerically.
\end{remk}

\noindent 
Before proving the theorem, we first apply it to our running example.
\begin{ex}
	In \Cref{ex:leadingreluctant}, we investigated a bifurcation scenario with a fully synchrony-breaking branch in the core, given by
	$x_i(\lambda) = D_i\lambda + \mathcal{O}(|\lambda|^2)$ for $i \in \{0,1,2\}$ and with mutually distinct $D_i$. 
	It follows that 
	$x_i(\lambda)- x_j(\lambda) = (D_i-D_j)\lambda + \mathcal{O}(|\lambda|^2)$ for all $i,j$. As $D_i-D_j\ne0$ for $i>j$, we obtain
	$p_{3,1} = p_{3,2} = p_{2,1} = 1$ and hence
	\begin{equation*}
		\bar{p} = 1+1+1 = 3\, .
	\end{equation*}
	\Cref{mainbifurcation} therefore predicts a bifurcation branch 
	\[ z(\lambda) = (x_0(\lambda), x_1(\lambda), x_2(\lambda), y_0(\lambda), y_1(\lambda)) \]
	in the augmented system $f^{\NNd}_{F,G}$, for a generic choice of $F$ and with any $G$ supporting the aforementioned fully synchrony-breaking branch in the core, which satisfies 
	\begin{equation*}
		y_0(\lambda) - y_1(\lambda) \sim \lambda^3\, .
	\end{equation*}
	This is indeed what we found in our numerical investigation, see \Cref{firstnum-c}. 
\end{ex}

\noindent The proof of \Cref{mainbifurcation} requires some machinery from \cite{vonderGracht.2023}. There we introduced the polynomials $P_{(k)}\colon \bigoplus_{\sigma \in S^0_{k+1}} \mathbb{R}^k \to  \mathbb{R}$, given by
\begin{equation}
	\label{eq:powersumsym1}
	P_{(k)}\left(\bigoplus_{\sigma \in S^0_{k+1}}  {\bf X}_{\sigma}\right) = \sum_{\sigma \in S^0_{k+1}}X_{\sigma,1}^{1}X_{\sigma,2}^{2}\dotsm X_{\sigma,k}^{k}\, ,
\end{equation}
for $k \in \N$, and where ${\bf{X}}_{\sigma} = (X_{\sigma,1}, \dots, X_{\sigma,k}) \in \R^k$ for $\sigma \in S^0_{k+1}$. 
We also state the following result, a proof of which can be found in \cite{vonderGracht.2023}.
\begin{lem}[Lemma 5.6 in \cite{vonderGracht.2023}]\label{lem:formain1}
	Let $Q\colon \bigoplus_{\sigma \in S_{k+1}^0} \R^k \rightarrow \R$ be a polynomial function that is invariant under all permutations of its $\#S_{k+1}^0$ entries from $\R^k$.
	Then there exists a polynomial $S \colon \R^{k+1} \rightarrow \R$ such that
	\begin{equation}\label{eq:factorisation}
		Q\left(\bigoplus_{\sigma \in S^0_{k+1}} {\bf x}_{\sigma} \right) - Q\left(\bigoplus_{\sigma \in S^1_{k+1}} {\bf x}_{\sigma} \right) = S(x) \prod_{\substack{i,j=0 \\ i > j }}^k (x_i - x_j)
	\end{equation}
	for all $x = (x_0, \dots, x_k) \in \R^{k+1}$, where ${\bf x}_{\sigma} = (x_{\sigma(1)}, \dots, x_{\sigma(k)})$ for all $\sigma \in S_{k+1}$.
\end{lem}

\begin{remk}\label{remk-samecompogivesamy}
	It can readily be seen that for any polynomial $Q$ satisfying the conditions of \Cref{lem:formain1}, \eqref{eq:factorisation} is actually equivalent to the fact that
	\begin{equation}\label{eq:similiartwothesame}
		Q\left(\bigoplus_{\sigma \in S^0_{k+1}} {\bf x}_{\sigma} \right) = Q\left(\bigoplus_{\sigma \in S^1_{k+1}} {\bf x}_{\sigma} \right)
	\end{equation}
	whenever $(x_0, \dots, x_k) \in \R^{k+1}$ satisfies $x_i = x_j$ for some distinct $i,j \in \{0,,\dots, k\}$. 
	This observation still holds when $Q$ is not polynomial, see Lemma 5.5 of \cite{vonderGracht.2023}. The latter fact actually  underlies the proof of \Cref{lem:formain1} that is given in \cite{vonderGracht.2023}.
\end{remk}

\begin{lem}\label{lem:formain2}
	The polynomials $P_{(k)}$ defined in \eqref{eq:powersumsym1} satisfy
	\begin{equation}\label{eq:factorisation_spec}
		P_{(k)}\left(\bigoplus_{\sigma \in S^0_{k+1}} {\bf x}_{\sigma} \right) - P_{(k)}\left(\bigoplus_{\sigma \in S^1_{k+1}} {\bf x}_{\sigma} \right) = \prod_{\substack{i,j=0 \\ i > j }}^k (x_i - x_j)
	\end{equation} 
	for all $x = (x_0, \dots, x_k) \in \R^{k+1}$.
\end{lem}

\begin{proof}
	By \Cref{lem:formain1} we have
	\begin{equation*}
		P_{(k)}\left(\bigoplus_{\sigma \in S^0_{k+1}} {\bf x}_{\sigma} \right) - P_{(k)}\left(\bigoplus_{\sigma \in S^1_{k+1}} {\bf x}_{\sigma} \right) = S(x)\prod_{\substack{i,j=0 \\ i > j }}^k (x_i - x_j)
	\end{equation*}
	for some polynomial $S$. 
	It remains to show that $S = 1$. 
	To this end, note that both the left and right hand side of \Cref{eq:factorisation_spec} has total degree $1+\dots+k = \frac{k(k+1)}{2}$. 
	This means $S$ is a constant polynomial. 
	As both sides of \Cref{eq:factorisation_spec} contain a term $1\cdot x_1x_2^2\dots x_k^k$, we see that $S = 1$ and the result follows.
\end{proof} 

\noindent Before we move on to the proof of \Cref{mainbifurcation}, we first have a closer look at the set of powers $\Upsilon$. Recall that we may write 
\begin{equation}
	\label{eq:formalpowerseries2}
	x_i(\lambda) = \sum_{p \in \Upsilon} A_{i,p} \lambda^{p} + \mathcal{O}(|\lambda|^{\bar{p}+\epsilon})\, 
\end{equation}
for all the components $x_i(\lambda)$ of a fully synchrony-breaking branch. 
By adding zero-co\-ef\-fi\-cients $A_{i,p}$ to Expression \eqref{eq:formalpowerseries2} and by decreasing $\epsilon$ if needed, we may assume that for all $p,q \in \Upsilon$, we have 
\begin{equation}
	\label{ex:suminvofupsilon}
	\begin{array}{rclclcl}
		p+q &\in& \Upsilon &\text{ if }& \, p+q &\leq& \bar{p}\,; \\ 
		p+q &\geq& \bar{p}+\epsilon  &\text{ if }& \, p+q &>& \bar{p}\,, 
	\end{array}
\end{equation}
and also that  
\begin{equation*}
	\begin{array}{rclclcl}
		1 &\in& \Upsilon &\text{ if }&\, 1 &\leq& \bar{p}\,; \\ 
		1 &\geq& \bar{p}+\epsilon  &\text{ if }&\, 1 &>& \bar{p}\,.
	\end{array}
\end{equation*}
More precisely, we can add to $\Upsilon$ all non-zero sums $s= c+\sum_{p \in \Upsilon}c_p p $ with non-negative integer coefficients $c,c_p$, such that $s \leq \bar{p}$. 
Note that this adds a finite number of elements to $\Upsilon$, as necessarily $c_p \leq \ceil*{\frac{\bar{p}}{p}}$ and  $c \leq \ceil*{\bar{p}}$. It then follows from \Cref{ex:suminvofupsilon} that $\bar{p} \in \Upsilon$, as we have $p_{i,j} \in \Upsilon$ for all $i>j$.
This allows us to iteratively investigate coefficients corresponding to (possibly non-integer) powers of $\lambda$ in the branches, as well as in polynomial expressions involving the components of these branches.
For instance, if $x_1(\lambda)$ and $x_2(\lambda)$ are given by \Cref{eq:formalpowerseries2}, then we may likewise write
\begin{align*}
	x_1(\lambda)x_2(\lambda) &= \sum_{p \in \Upsilon} A'_{p} \lambda^{p} + \mathcal{O}(|\lambda|^{\bar{p}+\epsilon}) \text{ and } \\ 
	\lambda x_1(\lambda) &= \sum_{p \in \Upsilon} A''_{p} \lambda^{p} + \mathcal{O}(|\lambda|^{\bar{p}+\epsilon})\, 
\end{align*}
for some $A'_{p}, A''_{p} \in \R$. 
Finally, whenever 
\begin{equation*}
	w(\lambda) = \sum_{p \in \Upsilon} A_{p} \lambda^{p} + \mathcal{O}(|\lambda|^{\bar{p}+\epsilon})\, 
\end{equation*}
for some locally defined map $w\colon \R_{\geq 0} \mapsto \R$ and with $A_p \in \R$, then for $q \in \Upsilon$ we may write 
\begin{equation*}
	[w(\lambda)]_{\leq q} := \sum_{\substack{p \in \Upsilon \\ p \leq q }} A_{p} \lambda^{p} \, \text{ and } \, 
	[w(\lambda)]_{< q} := \sum_{\substack{p \in \Upsilon \\ p < q }} A_{p}\lambda^{p}
\end{equation*}
for the truncated power series.

\begin{proof}[Proof of \Cref{mainbifurcation}]
	By assumption, $x(\lambda) = (x_0(\lambda), \dots,  x_k(\lambda))$ locally solves 
	\[ (f^{\NNd}_{(F,G)}(x,y;\lambda))_{v} = (f^{\NN}_{(G)}(x;\lambda))_v = 0\]
	for all nodes $v$ in the core
	$\NN$ and all $y = (y_0, y_1) \in \R^2$.
	To solve for the $y$-components, let $K\in \N$ be such that  
	\begin{equation}
		(K+1)\min(p \mid p \in \Upsilon) > \bar{p}\, .
	\end{equation}
	We expand a general response function $F$ as
	\begin{equation}\label{eq:iterativesolve0}
		F(Y,{\bf{X}};\lambda) = aY + \sum_{\ell,m=0}^K Q_{\ell,m}\left({\bf{X}}\right)Y^\ell\lambda^m + \mathcal{O}( \| \left({\bf{X}},Y;\lambda\right)\|^{K+1}) \,
	\end{equation}
	for $Y \in \R$, $\lambda \in \R_{\geq 0}$, and where
	\begin{equation*}
		{\bf{X}} := \bigoplus_{\sigma \in S^0_{k+1}} {\bf{X}}_{\sigma}
	\end{equation*}
	with ${\bf{X}}_{\sigma} \in \R^k$.
	Here each $Q_{\ell,m}$ is a polynomial of degree at most $K$ that is invariant under all permutations of the vectors ${\bf{X}}_{\sigma}$, which follows from the fact that $F$ is invariant under permutations of these vectors.  
	Our  assumption (which is necessary for a bifurcation) that $F(0,0;0) = 0$ implies that $Q_{0,0}(0) = 0$. 
	Moreover, by setting the number $a \in \R$ equal to the derivative of $F$ at $(0,0;0)$ in the $Y$-direction, we may assume that  $Q_{1,0}(0) = 0$. 
	
	For $s \in \{0,1\}$, the equation $\dot{y}_s = 0$ gives 
	\begin{equation}\label{eq:iterativesolve}
		ay_s + \sum_{\ell,m=0}^K Q_{\ell,m}\left(\bigoplus_{\sigma \in S^s_{k+1}} {\bf x}_{\sigma}(\lambda)\right)y_s^\ell\lambda^m + \mathcal{O}(\|(x(\lambda),y_s;\lambda)\|^{K+1}) = 0 \, ,
	\end{equation}
	where ${\bf x}_{\sigma}(\la) = (x_{\sigma(1)}(\la), \dots, x_{\sigma(k)}(\la))$ for all $\sigma \in S_{k+1}$.
	If we assume $a\not= 0$ then by the implicit function theorem, this equation locally has a unique solution $y_s(\lambda)$, which can be written as
	\begin{equation*}
		y_s(\lambda) = \sum_{p \in \Upsilon} B_{s,p} \lambda^{p} + \mathcal{O}(|\lambda|^{\bar{p}+\epsilon})\, .
	\end{equation*}
	Note that the lowest value of $p \in \Upsilon$ for which $B_{s,p} \not= 0$ will generically be the minimum of the orders of the $x_i(\lambda)$, unless these all exceed $1$. 
		In that case this lowest value of $p$ generically equals $1$, due to the presence of the term $Q_{0,1}(0)\lambda^1$ in Equation \eqref{eq:iterativesolve}. Thus by definition of $\hat{p}$ we indeed see that generically $y_s(\lambda) \sim \lambda^{\hat{p}}$.
	The coefficients $B_{s,p} \in \R$ can iteratively be solved for from the equation
	\begin{equation}\label{eq:iterativesolve2}
		ay_s(\lambda) + \sum_{\ell,m=0}^K Q_{\ell,m}\left(\bigoplus_{\sigma \in S^s_{k+1}} {\bf x}_{\sigma}(\lambda)\right)y_s^\ell(\lambda)\lambda^m + \mathcal{O}(|\lambda|^{\bar{p}+\epsilon}) = 0 \, ,
	\end{equation}
	which is \eqref{eq:iterativesolve} applied to the solution branch $(x_0(\la),\dotsc,x_k(\la),y_0(\la),y_1(\la))$.
	We now want to show that 
	\begin{equation}\label{eq:toshowdiffpowersyy}
		y_0(\lambda) - y_1(\lambda) = \mathcal{O}(|\lambda|^{\bar{p}})
	\end{equation}
	for these unique solutions. We will do so by proving for all $q \in \Upsilon$ with $q < \bar{p}$ that
	\begin{equation}\label{eq:toshowimplic}
		[y_0]_{<q} = [y_1]_{<q} \implies [y_0]_{\leq q} = [y_1]_{\leq q}\, .
	\end{equation}
	Note that for $q = \min(p \mid p \in \Upsilon)$ we have $[y_0]_{<q} = [y_1]_{<q} = 0$. 
	Hence, iterated use of Implication \eqref{eq:toshowimplic} indeed proves \Cref{eq:toshowdiffpowersyy}. 
	To show that the statement in \eqref{eq:toshowimplic} holds, we subtract \Cref{eq:iterativesolve2} for $s=1$ from the one for $s=0$, which gives us
	\begin{equation}
		a(y_0(\lambda)-y_1(\lambda)) + \sum_{\ell,m=0}^K \sum_{s = 0}^1 (-1)^s Q_{\ell,m}\left(\bigoplus_{\sigma \in S^s_{k+1}} {\bf x}_{\sigma}(\lambda)\right)y_s^\ell(\lambda)\lambda^m + \mathcal{O}(|\lambda|^{\bar{p}+\epsilon}) = 0 \, .
	\end{equation}
	Given $q \in \Upsilon$ satisfying $q < \bar{p}$, let $q^{+}$ denote the smallest element in $\Upsilon$ such that $q^{+} > q$, i.e. $q^{+}$ is the ``next power'' to consider. 
	Note that $q < \bar{p}$ means $q^{+} \leq \bar{p}$ exists. It follows that
	\begin{equation}\label{eq:iterativesolve4}
		\begin{split}
			0 &=   a([y_0(\lambda)]_{\leq q}-[y_1(\lambda)]_{\leq q}) \\  &\phantom{=}+ \sum_{\ell,m=0}^K \sum_{s = 0}^1 (-1)^s Q_{\ell,m}\left(\bigoplus_{\sigma \in S^s_{k+1}} {\bf x}_{\sigma}(\lambda)\right)[y_s(\lambda)]_{<q}^\ell\lambda^m + \mathcal{O}(|\lambda|^{q^+})  \, .
		\end{split}
	\end{equation}
	Here we have used that 
	\begin{multline}
		Q_{\ell,m}\left(\bigoplus_{\sigma \in S^s_{k+1}} {\bf x}_{\sigma}(\lambda)\right)[y_s(\lambda)]_{\leq q}^\ell\lambda^m \\
		= \, Q_{\ell,m}\left(\bigoplus_{\sigma \in S^s_{k+1}} {\bf x}_{\sigma}(\lambda)\right)[y_s(\lambda)]_{< q}^\ell\lambda^m + \mathcal{O}(|\lambda|^{q^+}) \, ,
	\end{multline}
	which is clear whenever $\ell = 0$, $\ell > 1$ or $m >0$. For $(\ell, m) = (1,0)$ it holds because $Q_{1,0}(0) = 0$, so that $Q_{1,0}$ has no constant term and hence
	$$Q_{1,0}\left(\bigoplus_{\sigma \in S^s_{k+1}}{\bf x}_{\sigma}(\lambda)\right)$$
	is divisible by $\lambda^r$ for $r = \min(p \mid p \in \Upsilon)$.
	We now assume $[y_0(\lambda)]_{< q} = [y_1(\lambda)]_{< q}$, so that $[y_s(\lambda)]_{< q} = [y_0(\lambda)]_{< q}$ for both choices of $s$. Using \Cref{lem:formain1}, \Cref{eq:iterativesolve4} becomes
	\begin{equation}\label{eq:iterativesolve5}
		\begin{split}
			0 &=   a([y_0(\lambda)]_{\leq q}-[y_1(\lambda)]_{\leq q}) \\  &+ \sum_{\ell,m=0}^K \sum_{s = 0}^1 (-1)^s Q_{\ell,m}\left(\bigoplus_{\sigma \in S^s_{k+1}} {\bf x}_{\sigma}(\lambda)\right)[y_s(\lambda)]_{<q}^\ell\lambda^m + \mathcal{O}(|\lambda|^{q^+}) \\ 
			&=   a([y_0(\lambda)]_{\leq q}-[y_1(\lambda)]_{\leq q}) \\  &+ \sum_{\ell,m=0}^K [y_0(\lambda)]_{<q}^\ell\lambda^m\sum_{s = 0}^1 (-1)^s Q_{\ell,m}\left(\bigoplus_{\sigma \in S^s_{k+1}} {\bf x}_{\sigma}(\lambda)\right) + \mathcal{O}(|\lambda|^{q^+}) \\ 
			&=   a([y_0(\lambda)]_{\leq q}-[y_1(\lambda)]_{\leq q}) \\  &+ \sum_{\ell,m=0}^K [y_0(\lambda)]_{<q}^\ell\lambda^m 
			S_{\ell, m}(x(\lambda)) \prod_{\substack{i,j=0 \\ i > j }}^k (x_i(\lambda) - x_j(\lambda)) + \mathcal{O}(|\lambda|^{q^+})\, , 
		\end{split}
	\end{equation}
	for some polynomials $S_{\ell, m}$. As it is clear that 
	\begin{equation*}
		\prod_{\substack{i,j=0 \\ i > j }}^k (x_i(\lambda) - x_j(\lambda)) = \mathcal{O}(|\lambda|^{\bar{p}})\, ,
	\end{equation*}
	\Cref{eq:iterativesolve5} simplifies to 
	\begin{equation*}
		a([y_0(\lambda)]_{\leq q}-[y_1(\lambda)]_{\leq q}) = \mathcal{O}(|\lambda|^{q^+})\, . 
	\end{equation*}
	Using again the assumption that $a\not= 0$, this indeed gives $[y_0(\lambda)]_{\leq q}=[y_1(\lambda)]_{\leq q}$. 
	
	By induction, \eqref{eq:toshowdiffpowersyy} holds true as outlined above.
	It follows that we may write
	\begin{equation*}
		y_0(\lambda) - y_1(\lambda) = E \lambda^{\bar{p}} + \mathcal{O}({|\lambda|^{\bar{p}+\epsilon}})\,  .
	\end{equation*} 
	for some $E= B_{0, \bar{p}} - B_{1,\bar{p}} \in \R$. 
	We next want to show that $E\not= 0$ generically. 
	To this end, recall that $B_{s,\bar{p}}$ can be solved for from \Cref{eq:iterativesolve2}. 
	In fact, for fixed values of $a\not= 0$ and the power series coefficients of each $x_i(\lambda)$, we may express  $B_{s,\bar{p}}$ as a polynomial in the coefficients of the various $Q_{\ell,m}$.
	Therefore, we may likewise express $E= B_{0, \bar{p}} - B_{1,\bar{p}}$ as such a polynomial.
	Now, any polynomial on a finite dimensional vector space is either identically zero, or vanishes only on the complement of an open dense set. 
	Therefore, the proof is complete if we can give at least one response function $F$ for which $E \not= 0$.
	To this end, consider 
	\begin{equation}
		F(Y,{\bf{X}};\lambda) = aY  + P_{(k)}({\bf{X}})\, .
	\end{equation}
	Using \Cref{lem:formain2} we get
	\begin{equation}
		\begin{split}
			0 &=  F\left(y_0(\la), \bigoplus_{\sigma \in S^0_{k+1}} {\bf x}_{\sigma}(\lambda); \lambda\right) - F\left(y_1(\la), \bigoplus_{\sigma \in S^1_{k+1}} {\bf x}_{\sigma}(\lambda); \lambda\right) \\ 
			&= a(y_0(\la) - y_1(\la)) +  P_{(k)}\left(\bigoplus_{\sigma \in S^0_{k+1}} {\bf x}_{\sigma}(\lambda) \right) - P_{(k)}\left(\bigoplus_{\sigma \in S^1_{k+1}} {\bf x}_{\sigma}(\lambda) \right) \\ 
			&= a(y_0(\la) - y_1(\la)) + \prod_{\substack{i,j=0 \\ i > j }}^k (x_i(\lambda) - x_j(\lambda)) \\ 
			&= a(y_0(\la) - y_1(\la)) + \bigg(\prod_{\substack{i,j=0 \\ i > j }}^k D_{i,j}\bigg)\lambda^{\bar{p}} + \mathcal{O}(|\lambda|^{\bar{p} + \epsilon})\, .
		\end{split}
	\end{equation}
	Hence, for this particular choice of response function we obtain 
	\begin{equation*}
		E = -\frac{1}{a}\prod_{\substack{i,j=0 \\ i > j }}^k D_{i,j} \not= 0\, ,
	\end{equation*}
	which shows that $E$ is indeed generically non-vanishing. 
	
	Finally, \Cref{eq:iterativesolve} shows that the branch 
	\begin{equation}
		z(\lambda) = (x_0(\lambda), \dots,  x_k(\lambda),  y_0(\lambda), y_1(\lambda))
	\end{equation}
	is stable if $x(\lambda)$ is stable for $f^{\NN}_{(G)}$, and if in addition $a<0$. Since $a=\frac{\partial F(0;0)}{\partial y}$, this completes the proof.
\end{proof}

\begin{remk}
	The condition in \Cref{mainbifurcation} that $x(\lambda)$ is fully synchrony-breaking is  essential. 
	More precisely, it follows from \Cref{remk-samecompogivesamy} that if $x_i(\lambda) = x_j(\lambda)$ for some distinct $i,j$, then the equations $\dot{y}_0 = 0$ and $\dot{y}_1 = 0$ give identical solutions $y_0(\lambda)$ and $y_1(\lambda)$. 
\end{remk}

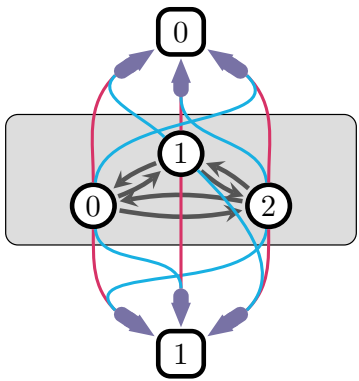
\begin{figure}
	\centering
	\resizebox{.35\linewidth}{!}{
		\begin{tikzpicture}[
	square/.style = {
		regular polygon,
		regular polygon sides=4
	},
	main node/.style={
		line width=1.5pt, 
		circle,
		draw,
		font=\sffamily,
		inner sep=2pt,
		fill=white
	},
	second node/.style={
		line width=1.5pt, 
		square, 
		draw, 
		font=\sffamily,
		rounded corners, 
		inner sep=2pt
	},
	edge/.style={
		-stealth,
		shorten >=1pt,
		shorten <=1pt,
		line width=1.5pt
	},
	hyperedge/.style={
		Round Cap-{Triangle[length=3mm, width=2mm]},
		line width=5pt
	}
	]
	\node[main node] (v0) at (-1,-.3) {$0$};
	\node[main node] (v1) at (0,.3) {$1$};
	\node[main node] (v2) at (1,-.3) {$2$};
	
	\node[second node] (w0) at (0,1.7) {$0$};
	\node[second node] (w1) at (0,-2) {$1$};
	
	\path[line width=1.5pt]

	(v1) edge [edge, grey, bend right=10] (v0)
	(v2) edge [edge, grey, bend right=10] (v0)
	
	(v0) edge [edge, grey, bend right=10] (v1)
	(v2) edge [edge, grey, bend right=10] (v1)
	
	(v0) edge [edge, grey, bend right=10] (v2)
	(v1) edge [edge, grey, bend right=10] (v2)
	;
	
	\coordinate (he1-1) at (-.75,1.2);
	\coordinate (he1-2) at (0,.95);
	\coordinate (he1-3) at (.75,1.2);
	\coordinate (he1-4) at (-.75,-1.5);
	\coordinate (he1-6) at (0,-1.25);
	\coordinate (he1-5) at (.75,-1.5);
	
	\path[line width=1pt]
	(v0) edge [red, in=south west, out=north] (he1-1)
	(v1) edge [blue, in=south west, out=north west] (he1-1)
	
	(v1) edge [red, in=south, out=north] (he1-2)
	(v2) edge [blue, in=south, out=95] (he1-2)
	
	(v2) edge [red, in=south east, out=north] (he1-3)
	(v0) edge [blue, in=south east, out=85] (he1-3)
	
	(v0) edge [red, in=north west, out=south] (he1-4)
	(v2) edge [blue, in=north west, out=265] (he1-4)
	
	(v2) edge [red, in=north east, out=south] (he1-5)
	(v1) edge [blue, in=north east, out=south east] (he1-5)
	
	(v1) edge [red, in=north, out=south] (he1-6)
	(v0) edge [blue, in=north, out=275] (he1-6)
	;

	\draw[hyperedge, mix] (he1-1) -- (w0);
	\draw[hyperedge, mix] (he1-2) -- (w0);
	\draw[hyperedge, mix] (he1-3) -- (w0);
	\draw[hyperedge, mix] (he1-4) -- (w1);
	\draw[hyperedge, mix] (he1-5) -- (w1);
	\draw[hyperedge, mix] (he1-6) -- (w1);
	
	\begin{pgfonlayer}{background}
		\draw[rounded corners, fill=grey!20] (-1.5,.75) -- (2,.75) -- (2,-.75) -- (-2,-.75) -- (-2,.75) -- (-1.5,.75);
	\end{pgfonlayer}
\end{tikzpicture}
	}
	\caption{An augmented hypernetwork that generically does not support reluctant synchrony-breaking.}
	\label{fig:counterex}
\end{figure}
\begin{ex}
	The augmented hypernetwork depicted in \Cref{fig:counterex} generically does not support a reluctant synchrony-breaking steady state branch. Its core (shown in the grey box) is a fully  symmetric $3$-cell network. It is known that this network only admits local synchrony breaking steady state branches for which $x_i(\lambda)=x_j(\lambda)$ for a pair $i\neq j$, see \cite{Golubitsky.1988}. Hence, all these branches are partially synchronous, and the conclusion of \Cref{mainbifurcation} cannot be drawn. According to the previous remark, all generic solution branches in fact satisfy $y_0(\la)=y_1(\la)$.
\end{ex}

\begin{remk}\label{remk:firstclosing}
	Recall from \Cref{moregeneralaug} that we may generalize the definition of an augmented hypernetwork to allow only hyperedges between the two additional nodes and a subset $\mathcal{S}$ of nodes of the same type of the core.
	It is clear that the results of \Cref{mainbifurcation} still hold for this construction.
	More precisely, we then need a steady-state bifurcation branch in the core that has different components for the nodes in $\mathcal{S}$. 
	We may then define the order of asynchrony as in \Cref{def:orderandorderofasyn}, but comparing only states in $\mathcal{S}$. 
	As in \Cref{mainbifurcation}, we will then generically have a steady-state bifurcation in a corresponding admissible vector field for the augmented hypernetwork, with the difference between the $w$-nodes growing in $\lambda$ raised to the power of the order of asynchrony.
	Stability of this branch for the augmented hypernetwork can be guaranteed if the relevant branch for the core is stable.
\end{remk}  

\section{More examples}
\label{sec:moreexamples}
\begin{figure}
	\begin{subfigure}{.35\linewidth}
		\centering
		\resizebox{.9\linewidth}{!}{
			\begin{tikzpicture}[
	square/.style = {
		regular polygon,
		regular polygon sides=4
	},
	main node/.style={
		line width=1.5pt, 
		circle,
		draw,
		font=\sffamily,
		inner sep=2pt,
		fill=white
	},
	second node/.style={
		line width=1.5pt, 
		square, 
		draw, 
		font=\sffamily,
		rounded corners, 
		inner sep=2pt
	},
	edge/.style={
		-stealth,
		shorten >=1pt,
		shorten <=1pt,
		line width=1.5pt
	},
	hyperedge/.style={
		Round Cap-{Triangle[length=3mm, width=2mm]},
		line width=5pt
	}
	]
	\node[main node] (v0) at (-1,0) {$0$};
	\node[main node] (v1) at (0,0) {$1$};
	\node[main node] (v2) at (1,0) {$2$};
	
	\node[second node] (w0) at (0,2) {$0$};
	\node[second node] (w1) at (0,-2) {$1$};
	
	\path[line width=1.5pt]

	(v1) edge [edge, grey, bend right=20] (v0)
	
	(v0) edge [edge, grey, bend right=20] (v1)
	
	(v2) edge [edge, grey, loop right, looseness=8] (v2)
	;
	
	\coordinate (he1-1) at (-.75,1.5);
	\coordinate (he1-2) at (0,1.1);
	\coordinate (he1-3) at (.75,1.5);
	\coordinate (he1-4) at (-.75,-1.5);
	\coordinate (he1-6) at (0,-1.1);
	\coordinate (he1-5) at (.75,-1.5);
	
	\path[line width=1pt]
	(v0) edge [red, in=south west, out=north] (he1-1)
	(v1) edge [blue, in=south west, out=north west] (he1-1)
	
	(v1) edge [red, in=south, out=north] (he1-2)
	(v2) edge [blue, in=south, out=north west] (he1-2)
	
	(v2) edge [red, in=south east, out=north] (he1-3)
	(v0) edge [blue, in=south east, out=north east] (he1-3)
	
	(v0) edge [red, in=north west, out=south] (he1-4)
	(v2) edge [blue, in=north west, out=south west] (he1-4)
	
	(v2) edge [red, in=north east, out=south] (he1-5)
	(v1) edge [blue, in=north east, out=south east] (he1-5)
	
	(v1) edge [red, in=north, out=south] (he1-6)
	(v0) edge [blue, in=north, out=south east] (he1-6)
	;

	\draw[hyperedge, mix] (he1-1) -- (w0);
	\draw[hyperedge, mix] (he1-2) -- (w0);
	\draw[hyperedge, mix] (he1-3) -- (w0);
	\draw[hyperedge, mix] (he1-4) -- (w1);
	\draw[hyperedge, mix] (he1-5) -- (w1);
	\draw[hyperedge, mix] (he1-6) -- (w1);
	
	\begin{pgfonlayer}{background}
		\draw[rounded corners, fill=grey!20] (-1,.75) -- (1.75,.75) -- (1.75,-.75) -- (-1.75,-.75) -- (-1.75,.75) -- (-1,.75);
	\end{pgfonlayer}
\end{tikzpicture}
		}
		\caption{An augmented hypernetwork with a disconnected core, shown within the grey box.}
		\label{fig:2nd_ex}
	\end{subfigure}
	\hfill
	\begin{subfigure}{.6\linewidth}
		\centering
		\includegraphics[width=\linewidth]{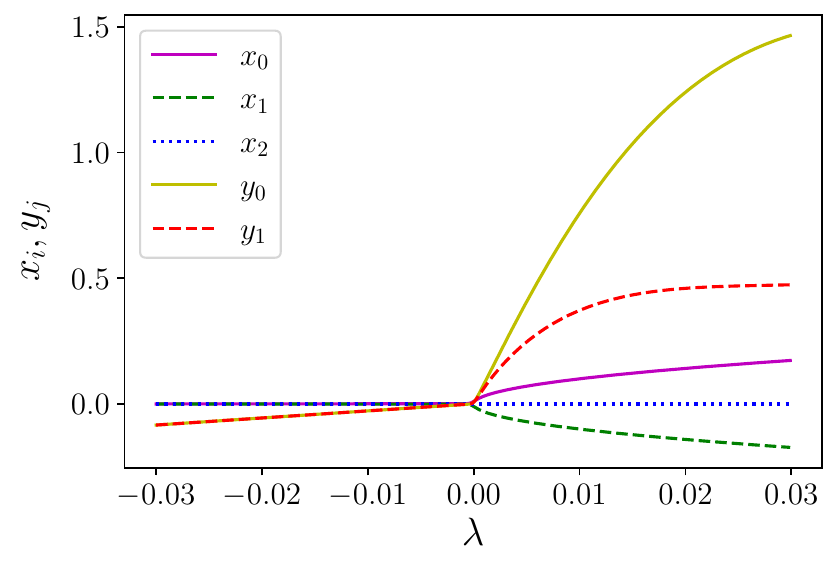}
		\caption{The stable branches of a synchrony-breaking bifurcation.}
		\label{secondnum-a}
	\end{subfigure}\par\medskip
	\begin{subfigure}{0.47\linewidth}
		\centering
		\includegraphics[width=\linewidth]{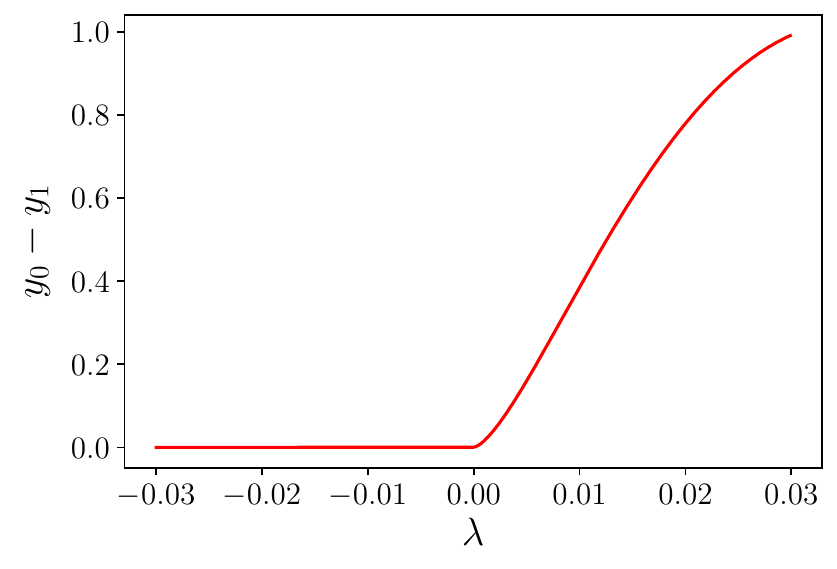}
		\caption{The difference between the $y$-nodes along the stable branches.}
		\label{secondnum-b}
	\end{subfigure}\hfill
	\begin{subfigure}{0.47\linewidth}
		\centering
		\includegraphics[width=\linewidth]{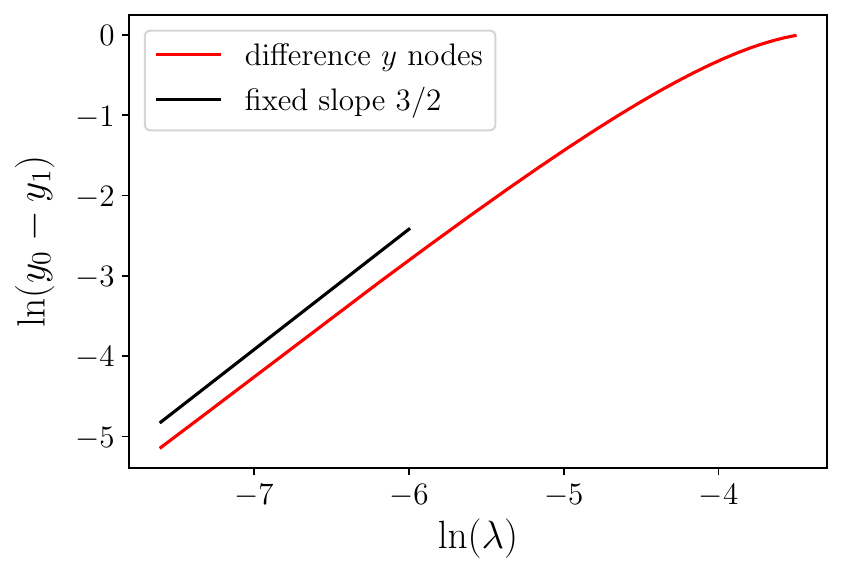}
		\caption{A log–log plot of the difference between the $y$-nodes.}
		\label{secondnum-c}
	\end{subfigure}
	\caption{An augmented hypernetwork with a disconnected core together with the numerically obtained bifurcation diagram for a corresponding system of the form \eqref{systembifwithpara2c}. We have left out self-loops corresponding to self-influence of each node. The black line segment in the log-log plot has fixed slope $\frac{3}{2}$ for comparison, indicating that $y_0(\lambda) - y_1(\lambda) \sim \lambda^{\frac{3}{2}}$. }
	\label{secondnum}
\end{figure}
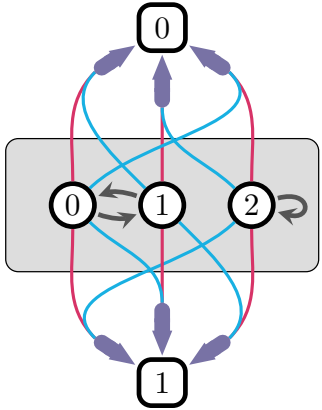
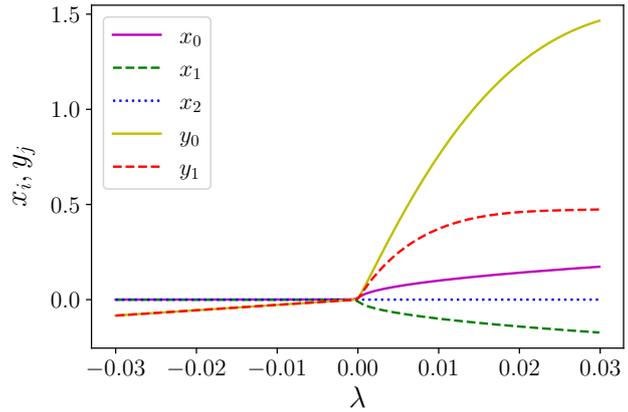
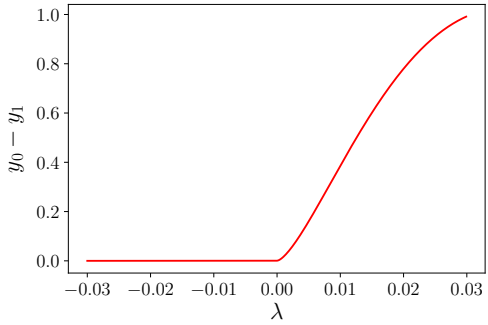
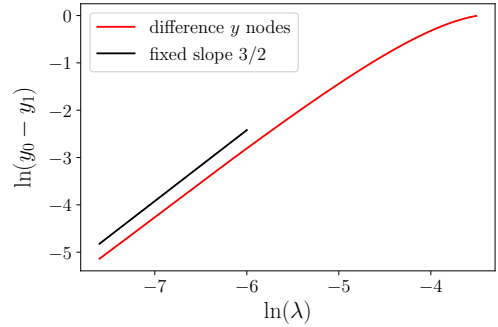

In this section, we present
three more examples that illustrate \Cref{mainbifurcation}.
\begin{ex}
	Consider the augmented hypernetwork shown in \Cref{fig:2nd_ex}. Its admissible ODEs are given by
	\begin{equation}\label{systembifwithpara2c}
		\begin{array}{ll}
			\dot x_0 & = G(x_0, \textcolor{grey}{x_1}; \lambda)\, , \\
			\dot x_1 & = G(x_1, \textcolor{grey}{x_0}; \lambda)\, , \\
			\dot x_2 & = G(x_2, \textcolor{grey}{x_2}; \lambda)\, , \\
			\dot y_0 & = F(y_0, \textcolor{mix}{(} \textcolor{red}{x_0}, \textcolor{blue}{x_1} \textcolor{mix}{)}, \textcolor{mix}{(} \textcolor{red}{x_1}, \textcolor{blue}{x_2} \textcolor{mix}{)}, \textcolor{mix}{(} \textcolor{red}{x_2}, \textcolor{blue}{x_0} \textcolor{mix}{)}; \lambda)\, ,\\
			\dot y_1 & = F(y_1, \textcolor{mix}{(} \textcolor{red}{x_0}, \textcolor{blue}{x_2} \textcolor{mix}{)}, \textcolor{mix}{(} \textcolor{red}{x_1}, \textcolor{blue}{x_0} \textcolor{mix}{)}, \textcolor{mix}{(} \textcolor{red}{x_2}, \textcolor{blue}{x_1} \textcolor{mix}{)}; \lambda)\, ,\\
		\end{array}
	\end{equation}
	where $F$ has the usual symmetry properties and where we assume each node to have a one-dimensional phase space. 
	The grey box in \Cref{fig:2nd_ex} denotes the core of this hypernetwork, which is a disconnected, classical first-order network, and whose dynamics corresponds to that of the $x$-variables in \eqref{systembifwithpara2c}.
	Generically, a one-parameter bifurcation in an admissible system for the core is either given by the product of two saddle-nodes or by a pitchfork bifurcation. 
	Only the latter of these involves a fully synchrony-breaking branch, and so we focus on that case, which we realise by choosing
	\begin{equation}\label{eq:Gin2ndex}
		G(X_0,X_1;\lambda) =  -X_0-X_1+\lambda X_0-X_0^3 \, .
	\end{equation}
	It follows that for $\lambda<0$, the only (stable) steady state branch is given by $x_0(\lambda) = x_1(\lambda) = x_2(\lambda) = 0$. For $\lambda>0$ we find (apart from some unstable branches) two stable, fully synchrony-breaking branches, given by
	\begin{equation*}
		x_0(\lambda) = - x_1(\lambda) = \pm \lambda^{1/2}, \, x_2(\lambda) = 0\, .
	\end{equation*}
	For each of these two latter branches we have 
	\begin{equation*}
		x_2(\lambda) - x_0(\lambda) =  \mp \lambda^{1/2}, \, \,
		x_2(\lambda) - x_1(\lambda) =  \pm \lambda^{1/2}, x_1(\lambda) - x_0(\lambda) =  \mp 2\lambda^{1/2}\, ,
	\end{equation*}
	from which we see that
	\begin{equation*}
		\bar{p} = \frac{1}{2}+\frac{1}{2}+\frac{1}{2} = \frac{3}{2}\, .
	\end{equation*}
	Hence, for the choice of response function $G$ given by \eqref{eq:Gin2ndex}, \Cref{mainbifurcation} predicts the full system \eqref{systembifwithpara2c} to undergo a bifurcation where $y_0(\lambda) - y_1(\lambda) \sim \lambda^{3/2}$ for $\lambda > 0$, for a generic choice of $F$. 
	
	A numerical investigation corroborates this result, see \Cref{secondnum}. These figures are obtained using Euler's method for the system \eqref{systembifwithpara2c} with $G$ given by \Cref{eq:Gin2ndex} and where we use
	\begin{align}\label{num:res-F}
		&F(Y, (X_0, X_1), (X_2, X_3), (X_4, X_5); {\lambda} ) \\ \nonumber =  &-5Y + 14 \lambda 
		-h(10X_0-12X_1)-h(10X_2-12X_3)-h(10X_4-12X_5)\,
	\end{align}
	in which
	\begin{equation}\label{num:res-h} 
		h(x) = \sin(x)+\cos(x)-1 = \sqrt{2}\sin(x+\frac{\pi}{4}) - 1\, .
	\end{equation}
	This function $F$ is chosen because it has non-vanishing Taylor coefficients of arbitrary order---to guarantee that the genericity conditions of \Cref{mainbifurcation} hold---, while also satisfying the required symmetry condition.
	We forward integrated the system \eqref{systembifwithpara2c} for each of $600$ equidistributed values of $\lambda \in [-0.03, 0.03]$. 
	For each fixed value of $\lambda$, integration was performed up to $t=5000$ with time steps of $0.1$, and starting from the point $(x_0, x_1, x_2, y_0, y_1) = (0.1, -0.2, 0.3, 0.4, 0.5)$.
	For the log-log plot of \Cref{secondnum-c}, we instead chose $600$ values of $\lambda \in [0.0005, 0.03]$, such that the values of $\ln(\lambda)$ are equidistributed.
\end{ex}

\begin{ex}
	\label{example:largefund}
	We next turn to the augmented hypernetwork depicted in \Cref{fig:3nd_ex}. The core of this hypernetwork, shown in the grey box, is an example of a classical (dyadic) network that itself shows reluctant synchrony-breaking. More precisely, admissible systems for this core are of the form
	\begin{equation}\label{system5fund}
		\begin{array}{ll}
			\dot x_0 & = G(x_0, \textcolor{uhhdarkgrey}{x_1}, \textcolor{grey}{x_0}; \lambda)\, , \\
			\dot x_1 & = G(x_1, \textcolor{uhhdarkgrey}{x_2}, \textcolor{grey}{x_0}; \lambda) \, , \\
			\dot x_2 & = G(x_2, \textcolor{uhhdarkgrey}{x_2}, \textcolor{grey}{x_0}; \lambda)\, . \\
		\end{array}
	\end{equation}
	These ODEs are special instances of those of the more general form 
	\begin{equation}\label{system5fund2}
		\begin{array}{ll}
			\dot x_0 & = H(x_0, \textcolor{uhhdarkgrey}{x_1}, \textcolor{grey}{x_0}, \textcolor{weirdgrey}{x_1}, \textcolor{notevengrey}{x_2}; \lambda)\, ,  \\
			\dot x_1 & = H(x_1, \textcolor{uhhdarkgrey}{x_2}, \textcolor{grey}{x_0}, \textcolor{weirdgrey}{x_1}, \textcolor{notevengrey}{x_2}; \lambda) \, , \\
			\dot x_2 & = H(x_2, \textcolor{uhhdarkgrey}{x_2}, \textcolor{grey}{x_0}, \textcolor{weirdgrey}{x_1}, \textcolor{notevengrey}{x_2}; \lambda)\, , \\
		\end{array}
	\end{equation}
	obtained by setting $H(x,y,z,u,v; \lambda) = G(x,y,z; \lambda)$. 
	Alternatively, one may think of \eqref{system5fund2} as denoting all admissible systems for a network obtained from the core in \Cref{fig:3nd_ex} by adding six additional arrows. 
	The first three of these are from node 1 to all nodes in the core (including an additional self-loop for node 1), and are all of a single, new type. 
	The last three are from node 2 to all nodes in the core, likewise all of one new type.
	The reason for adding these new arrow-types is that we can rigorously compute generic steady state bifurcations in systems of the form \eqref{system5fund2}, using center manifold reduction. See \cite{nijholt2019center} for a detailed exposition of the techniques used. 
	
	For the sake of this example, it is enough to know that in  Subsection 7.2 of \cite{nijholt2019center} it is shown that the system \eqref{system5fund2} generically undergoes a steady-state bifurcation involving a synchrony-breaking branch
	\begin{equation}\label{branches5cnetwork}
		\begin{split}
			x(\lambda) &= (x_0(\lambda), x_1(\lambda), x_2(\lambda)) \\ 
			&= (D_0\lambda + \mathcal{O}(|\lambda|^2), \, D_1\lambda + \mathcal{O}(|\lambda|^2), \, D_1\lambda + \mathcal{O}(|\lambda|^2))\, ,
		\end{split}
	\end{equation}
	where  
	\begin{equation}\label{branches5cnetwork-2}
		x_2(\lambda) - x_1(\lambda) = D_{2,1}\lambda^2 + \mathcal{O}(|\lambda|^3))\, 
	\end{equation}
	and with $D_0, D_1, D_0-D_1, D_{2,1} \not= 0$. 
	As this branch diverges from the synchrony space $\{x_1 = x_2\}$ at only quadratic leading order, we may again speak of reluctant synchrony breaking. 
	As opposed to the reluctant synchrony breaking we have considered in augmented hypernetworks though, the space $\{x_1 = x_2\}$ is actually robust for systems of the form \eqref{system5fund2}, and so for the special cases \eqref{system5fund} as well. 
	The synchrony-breaking branch \eqref{branches5cnetwork} can furthermore take over stability from a fully synchronous one as $\lambda$ increases through zero, see Table 2.1 in \cite{nijholt2019center}. 
	We therefore predict such a bifurcation to occur in the special system \eqref{system5fund} as well.
	
	\begin{figure}
		\begin{subfigure}{.35\linewidth}
			\centering
			\resizebox{.9\linewidth}{!}{
				\begin{tikzpicture}[
	square/.style = {
		regular polygon,
		regular polygon sides=4
	},
	main node/.style={
		line width=1.5pt, 
		circle,
		draw,
		font=\sffamily,
		inner sep=2pt,
		fill=white
	},
	second node/.style={
		line width=1.5pt, 
		square, 
		draw, 
		font=\sffamily,
		rounded corners, 
		inner sep=2pt
	},
	edge/.style={
		-stealth,
		shorten >=1pt,
		shorten <=1pt,
		line width=1.5pt
	},
	hyperedge/.style={
		Round Cap-{Triangle[length=3mm, width=2mm]},
		line width=5pt
	}
	]
	\node[main node] (v0) at (-1,0) {$0$};
	\node[main node] (v1) at (0,0) {$1$};
	\node[main node] (v2) at (1,0) {$2$};
	
	\node[second node] (w0) at (0,2) {$0$};
	\node[second node] (w1) at (0,-2) {$1$};
	
	\path[line width=1.5pt]

	(v1) edge [edge, uhhdarkgrey, bend right=20] (v0)
	
	(v2) edge [edge, uhhdarkgrey, bend right=20] (v1)
	
	(v2) edge [edge, uhhdarkgrey, loop right, looseness=8] (v2)

        (v0) edge [edge, grey, loop left, looseness=8] (v0)
	
	(v0) edge [edge, grey, bend right=16] (v1)
	
	(v0) edge [edge, grey, bend right=27] (v2)

	;
	
	\coordinate (he1-1) at (-.75,1.5);
	\coordinate (he1-2) at (0,1.1);
	\coordinate (he1-3) at (.75,1.5);
	\coordinate (he1-4) at (-.75,-1.5);
	\coordinate (he1-6) at (0,-1.1);
	\coordinate (he1-5) at (.75,-1.5);
	
	\path[line width=1pt]
	(v0) edge [red, in=south west, out=north] (he1-1)
	(v1) edge [blue, in=south west, out=north west] (he1-1)
	
	(v1) edge [red, in=south, out=north] (he1-2)
	(v2) edge [blue, in=south, out=north west] (he1-2)
	
	(v2) edge [red, in=south east, out=north] (he1-3)
	(v0) edge [blue, in=south east, out=north east] (he1-3)
	
	(v0) edge [red, in=north west, out=south] (he1-4)
	(v2) edge [blue, in=north west, out=south west] (he1-4)
	
	(v2) edge [red, in=north east, out=south] (he1-5)
	(v1) edge [blue, in=north east, out=south east] (he1-5)
	
	(v1) edge [red, in=north, out=south] (he1-6)
	(v0) edge [blue, in=north, out=south east] (he1-6)
	;

	\draw[hyperedge, mix] (he1-1) -- (w0);
	\draw[hyperedge, mix] (he1-2) -- (w0);
	\draw[hyperedge, mix] (he1-3) -- (w0);
	\draw[hyperedge, mix] (he1-4) -- (w1);
	\draw[hyperedge, mix] (he1-5) -- (w1);
	\draw[hyperedge, mix] (he1-6) -- (w1);
	
	\begin{pgfonlayer}{background}
		\draw[rounded corners, fill=grey!20] (-1,.75) -- (1.75,.75) -- (1.75,-.75) -- (-1.75,-.75) -- (-1.75,.75) -- (-1,.75);
	\end{pgfonlayer}
\end{tikzpicture}
			}
			\caption{An augmented hypernetwork with a core that itself shows reluctant synchrony breaking.}
			\label{fig:3nd_ex}
		\end{subfigure}
		\hfill
		\begin{subfigure}{.6\linewidth}
			\centering
			\includegraphics[width=\linewidth]{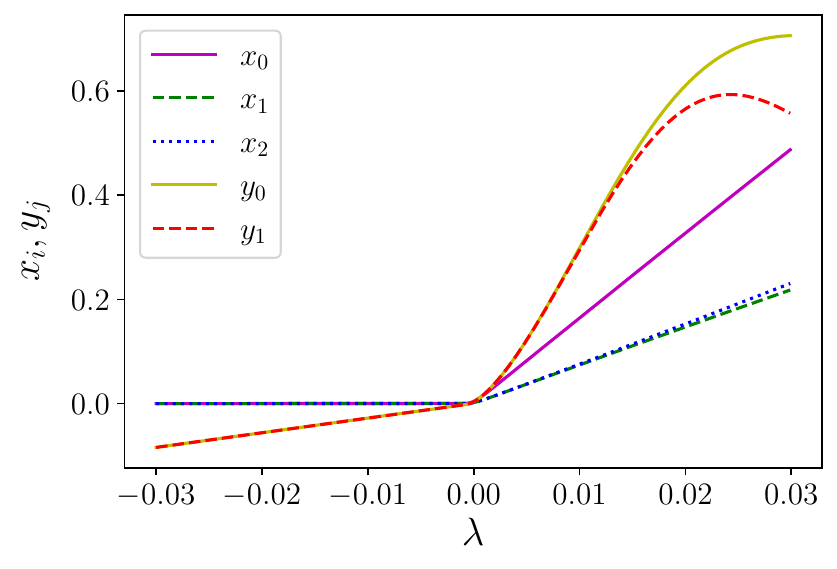}
			\caption{The stable branches of a synchrony-breaking bifurcation.}
			\label{third-a}
		\end{subfigure}\par\medskip
		\begin{subfigure}{0.47\linewidth}
			\centering
			\includegraphics[width=\linewidth]{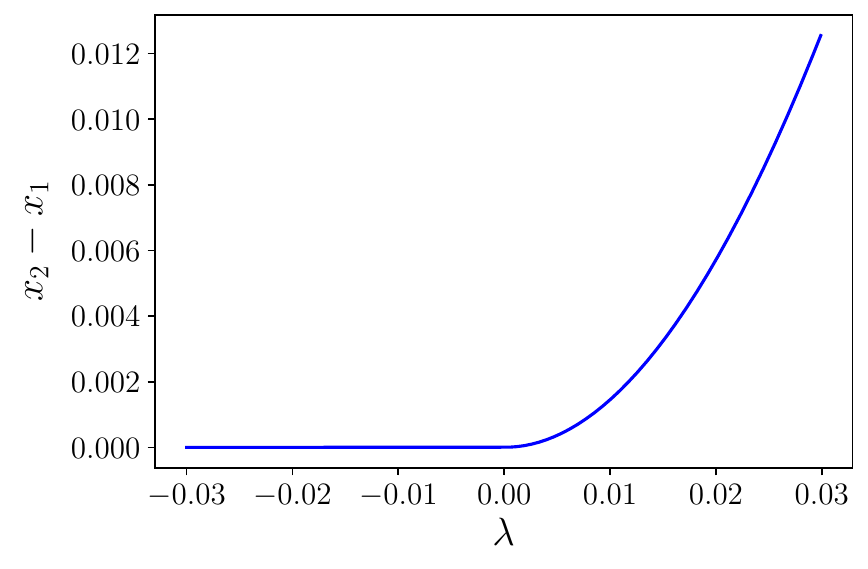}
			\caption{The difference between $x_2$ and $x_1$ along the stable branches.}
			\label{third-b}
		\end{subfigure}\hfill
		\begin{subfigure}{0.47\linewidth}
			\centering
			\includegraphics[width=\linewidth]{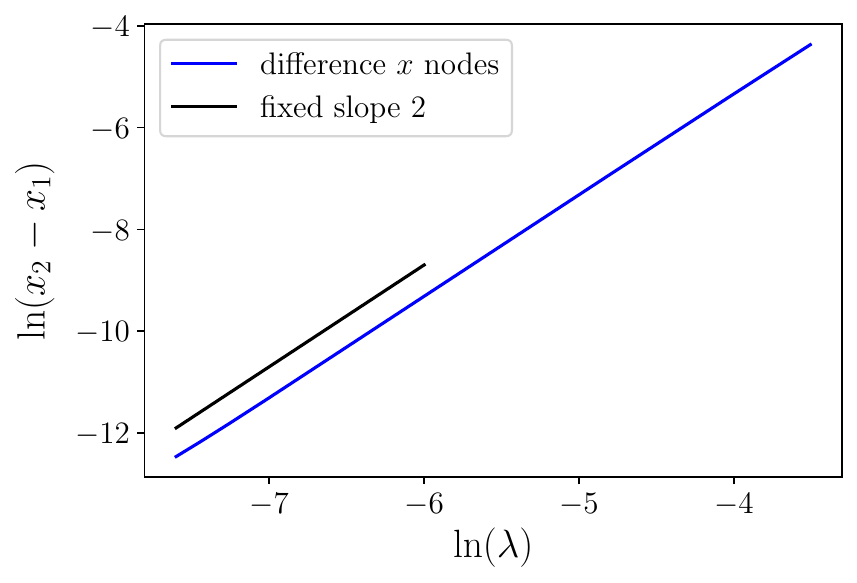}
			\caption{A log–log plot of the difference between between $x_2$ and $x_1$.}
			\label{third-c}
		\end{subfigure}
		\begin{subfigure}{0.47\linewidth}
			\centering
			\includegraphics[width=\linewidth]{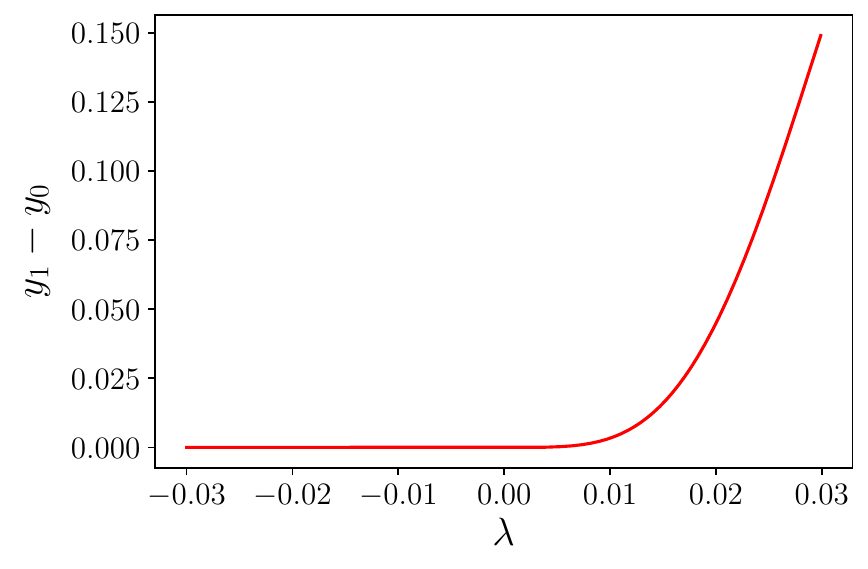}
			\caption{The difference between the  $y$-nodes along the stable branches.}
			\label{third-d}
		\end{subfigure}\hfill
		\begin{subfigure}{0.47\linewidth}
			\centering
			\includegraphics[width=\linewidth]{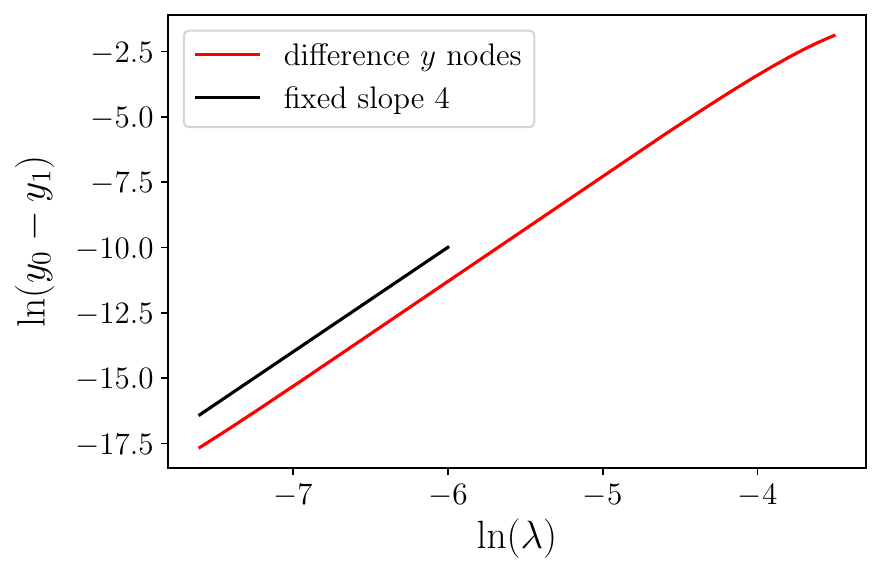}
			\caption{A log–log plot of the difference between the $y$-nodes.}
			\label{third-e}
		\end{subfigure}
		\caption{An augmented hypernetwork together with the numerically obtained bifurcation diagram whose corresponding core system is of the form \eqref{system5fund}. We have left out self-loops corresponding to self-influence of each node. The black line segment in (d) has fixed slope~$2$ for comparison, indicating that $x_2(\lambda) - x_1(\lambda) \sim \lambda^{2}$. The black line segment in (f) has fixed slope~$4$, which indicates that $y_0(\lambda) - y_1(\lambda) \sim \lambda^{4}$.}
		\label{thirdnum}
	\end{figure}
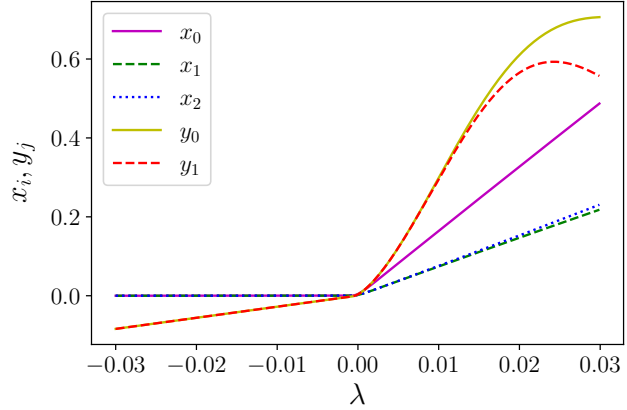
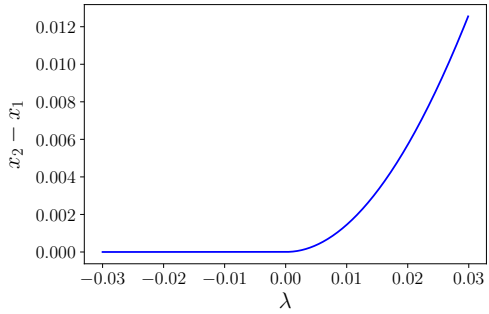
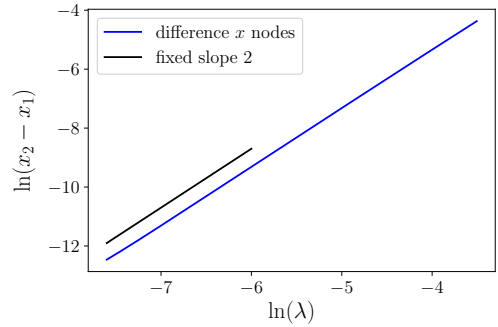
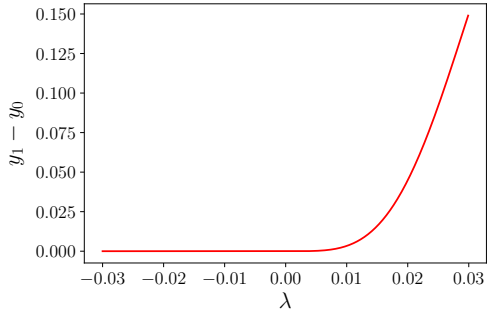
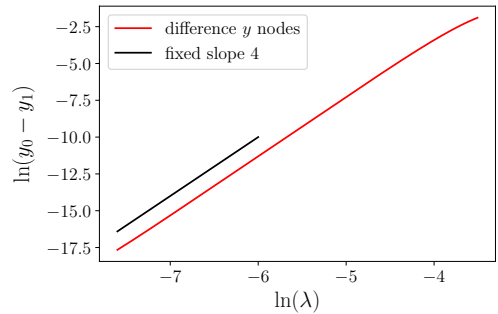
	
	\noindent \Cref{thirdnum} reveals that this is indeed the case. 
	\Cref{third-a} shows the components of the stable branches for the augmented hypernetwork of \Cref{fig:3nd_ex}. 
	However, as the core is a subnetwork of its augmented hypernetwork, we see that the $x$-variables depend only on each other and so depict a bifurcation in the three-node system \eqref{system5fund} as well.
	\Cref{third-a} shows a reluctant separation happening between the nodes $v_1$ and $v_2$, corresponding to the variables $x_1$ and $x_2$, with \Cref{third-b,third-c} indicating this occurs as $\sim \lambda^2$. 
	
	It follows from \eqref{branches5cnetwork} and \eqref{branches5cnetwork-2} that we have
	\begin{equation*}
		x_2(\lambda)-x_0(\lambda) \sim \la, \, x_1(\lambda)-x_0(\lambda) \sim \la \text{ and } x_2(\lambda)-x_1(\lambda) \sim \la^2\, 
	\end{equation*}
	so that 
	\begin{equation*}
		\bar{p} = 1 + 1 + 2 = 4\, .
	\end{equation*}
	By \Cref{mainbifurcation} this implies that the augmented hypernetwork system generically exhibits the highly reluctant synchrony-breaking asymptotics  
	\begin{equation*}
		y_0(\lambda) - y_1(\lambda) \sim \lambda^4\, .
	\end{equation*}
	The numerics in \Cref{third-d,third-e} corroborates this surprising asymptotics -- see in particular the branches corresponding to the $y$-variables in \Cref{third-a}.
	
	The details of the numerics are the same as for the previous example, except that time ran until $t=5000$ for \Cref{third-a,third-b,third-d}, and until $t=15000$ for \Cref{third-c,third-e}. The response function for the $x$-variables was chosen to be
	\begin{equation*}
		G(X_0,X_1,X_2;\lambda) =  -0.55X_1+0.25X_2+1.5\lambda X_0-0.1X_0^2 \, ,
	\end{equation*}
	with $F$ given by \Cref{num:res-F,num:res-h}.
\end{ex} 

\begin{figure}
	\begin{minipage}{.35\linewidth}
		\begin{subfigure}{\linewidth}
			\centering
			\resizebox{\linewidth}{!}{
				\begin{tikzpicture}[
	square/.style = {
		regular polygon,
		regular polygon sides=4
	},
	main node/.style={
		line width=1.5pt, 
		circle,
		draw,
		font=\sffamily,
		inner sep=2pt,
		fill=white
	},
	second node/.style={
		line width=1.5pt, 
		square, 
		draw, 
		font=\sffamily,
		rounded corners, 
		inner sep=2pt
	},
	edge/.style={
		-stealth,
		shorten >=1pt,
		shorten <=1pt,
		line width=1.5pt
	},
	hyperedge/.style={
		Round Cap-{Triangle[length=3mm, width=2mm]},
		line width=5pt
	},
	pol/.style = {
		regular polygon,
		regular polygon sides=36
	},
	pol node/.style={
		pol,
		draw=none,
		minimum size=2.5cm
	}
	]
	
	\clip(-2,-3.5) rectangle (2.3,3.5);
	
	\node[main node] (v0) at (-1.5,0) {$0$};
	\node[main node] (v1) at (-.5,0) {$1$};
	\node[main node] (v2) at (.5,0) {$2$};
	\node[main node] (v3) at (1.5,0) {$3$};
	
	\node[second node] (w0) at (0,3) {$0$};
	\node[second node] (w1) at (0,-3) {$1$};
	
	\node[pol node, rotate=180] (p0) at (0,3) {};	
	\coordinate (he0-1) at (p0.corner 32);
	\coordinate (he0-2) at (p0.corner 33);
	\coordinate (he0-3) at (p0.corner 34);
	\coordinate (he0-4) at (p0.corner 35);
	\coordinate (he0-5) at (p0.corner 36);
	\coordinate (he0-6) at (p0.corner 1);
	\coordinate (he0-7) at (p0.corner 2);
	\coordinate (he0-8) at (p0.corner 3);
	\coordinate (he0-9) at (p0.corner 4);
	\coordinate (he0-10) at (p0.corner 5);
	\coordinate (he0-11) at (p0.corner 6);
	\coordinate (he0-12) at (p0.corner 7);

	\node[pol node] (p1) at (0,-3) {};
	\coordinate (he1-1) at (p1.corner 7);
	\coordinate (he1-2) at (p1.corner 6);
	\coordinate (he1-3) at (p1.corner 5);
	\coordinate (he1-4) at (p1.corner 4);
	\coordinate (he1-5) at (p1.corner 3);
	\coordinate (he1-6) at (p1.corner 2);
	\coordinate (he1-7) at (p1.corner 1);
	\coordinate (he1-8) at (p1.corner 36);
	\coordinate (he1-9) at (p1.corner 35);
	\coordinate (he1-10) at (p1.corner 34);
	\coordinate (he1-11) at (p1.corner 33);
	\coordinate (he1-12) at (p1.corner 32);
	
	\path[line width=1pt]
	(v1) edge [red, in=south west, out=135] (he0-1)
	(v2) edge [blue, in=south west, out=135] (he0-1)
	(v3) edge [upborange, in=south west, out=135] (he0-1)
	
	(v0) edge [red, in=south west, out=90] (he0-2)
	(v3) edge [blue, in=south west, out=130] (he0-2)
	(v2) edge [upborange, in=south west, out=125] (he0-2)
	
	(v3) edge [red, in=south west, out=125] (he0-3)
	(v0) edge [blue, in=south west, out=85] (he0-3)
	(v1) edge [upborange, in=south west, out=125] (he0-3)
	
	(v2) edge [red, in=south west, out=115] (he0-4)
	(v1) edge [blue, in=south west, out=115] (he0-4)
	(v0) edge [upborange, in=south west, out=80] (he0-4)
	
	(v0) edge [red, in=south, out=75] (he0-5)
	(v1) edge [blue, in=south, out=105] (he0-5)
	(v3) edge [upborange, in=south, out=120] (he0-5)
	
	(v2) edge [red, in=south, out=105] (he0-6)
	(v0) edge [blue, in=south, out=70] (he0-6)
	(v3) edge [upborange, in=south, out=115] (he0-6)
	
	(v0) edge [red, in=south, out=65] (he0-7)
	(v2) edge [blue, in=south, out=95] (he0-7)
	(v1) edge [upborange, in=south, out=95] (he0-7)
	
	(v3) edge [red, in=south, out=110] (he0-8)
	(v2) edge [blue, in=south, out=85] (he0-8)
	(v0) edge [upborange, in=south, out=60] (he0-8)
	
	(v1) edge [red, in=south east, out=85] (he0-9)
	(v0) edge [blue, in=south east, out=55] (he0-9)
	(v2) edge [upborange, in=south east, out=75] (he0-9)
	
	(v1) edge [red, in=south east, out=75] (he0-10)
	(v3) edge [blue, in=south east, out=105] (he0-10)
	(v0) edge [upborange, in=south east, out=50] (he0-10)
	
	(v3) edge [red, in=south east, out=100] (he0-11)
	(v1) edge [blue, in=south east, out=65] (he0-11)
	(v2) edge [upborange, in=south east, out=65] (he0-11)
	
	(v2) edge [red, in=south east, out=55] (he0-12)
	(v3) edge [blue, in=south east, out=90] (he0-12)
	(v1) edge [upborange, in=south east, out=55] (he0-12)
	;

	\path[line width=1pt]
	(v0) edge [red, in=north west, out=-90] (he1-1)
	(v2) edge [blue, in=north west, out=-135] (he1-1)
	(v3) edge [upborange, in=north west, out=-135] (he1-1)
	
	(v1) edge [red, in=north west, out=-135] (he1-2)
	(v0) edge [blue, in=north west, out=-85] (he1-2)
	(v3) edge [upborange, in=north west, out=-125] (he1-2)
	
	(v1) edge [red, in=north west, out=-125] (he1-3)
	(v2) edge [blue, in=north west, out=-125] (he1-3)
	(v0) edge [upborange, in=north west, out=-75] (he1-3)
	
	(v2) edge [red, in=north west, out=-115] (he1-4)
	(v1) edge [blue, in=north west, out=-115] (he1-4)
	(v3) edge [upborange, in=north west, out=-120] (he1-4)
	
	(v3) edge [red, in=north, out=-115] (he1-5)
	(v2) edge [blue, in=north, out=-105] (he1-5)
	(v1) edge [upborange, in=north, out=-105] (he1-5)
	
	(v1) edge [red, in=north, out=-95] (he1-6)
	(v3) edge [blue, in=north, out=-110] (he1-6)
	(v2) edge [upborange, in=north, out=-95] (he1-6)
	
	(v0) edge [red, in=north, out=-70] (he1-7)
	(v1) edge [blue, in=north, out=-85] (he1-7)
	(v2) edge [upborange, in=north, out=-85] (he1-7)
	
	(v3) edge [red, in=north, out=-105] (he1-8)
	(v0) edge [blue, in=north, out=-65] (he1-8)
	(v2) edge [upborange, in=north, out=-75] (he1-8)
	
	(v3) edge [red, in=north east, out=-100] (he1-9)
	(v1) edge [blue, in=north east, out=-75] (he1-9)
	(v0) edge [upborange, in=north east, out=-60] (he1-9)
	
	(v0) edge [red, in=north east, out=-55] (he1-10)
	(v3) edge [blue, in=north east, out=-95] (he1-10)
	(v1) edge [upborange, in=north east, out=-65] (he1-10)
	
	(v2) edge [red, in=north east, out=-65] (he1-11)
	(v0) edge [blue, in=north east, out=-50] (he1-11)
	(v1) edge [upborange, in=north east, out=-55] (he1-11)
	
	(v0) edge [red, in=north east, out=-45] (he1-12)
	(v1) edge [blue, in=north east, out=-45] (he1-12)
	(v2) edge [upborange, in=north east, out=-55] (he1-12)
;

	\draw[hyperedge, newmix] (he0-1) -- (w0);
	\draw[hyperedge, newmix] (he0-2) -- (w0);
	\draw[hyperedge, newmix] (he0-3) -- (w0);
	\draw[hyperedge, newmix] (he0-4) -- (w0);
	\draw[hyperedge, newmix] (he0-5) -- (w0);
	\draw[hyperedge, newmix] (he0-6) -- (w0);
	\draw[hyperedge, newmix] (he0-7) -- (w0);
	\draw[hyperedge, newmix] (he0-8) -- (w0);
	\draw[hyperedge, newmix] (he0-9) -- (w0);
	\draw[hyperedge, newmix] (he0-10) -- (w0);
	\draw[hyperedge, newmix] (he0-11) -- (w0);
	\draw[hyperedge, newmix] (he0-12) -- (w0);
	
	\draw[hyperedge, newmix] (he1-1) -- (w1);
	\draw[hyperedge, newmix] (he1-2) -- (w1);
	\draw[hyperedge, newmix] (he1-3) -- (w1);
	\draw[hyperedge, newmix] (he1-4) -- (w1);
	\draw[hyperedge, newmix] (he1-5) -- (w1);
	\draw[hyperedge, newmix] (he1-6) -- (w1);
	\draw[hyperedge, newmix] (he1-7) -- (w1);
	\draw[hyperedge, newmix] (he1-8) -- (w1);
	\draw[hyperedge, newmix] (he1-9) -- (w1);
	\draw[hyperedge, newmix] (he1-10) -- (w1);
	\draw[hyperedge, newmix] (he1-11) -- (w1);
	\draw[hyperedge, newmix] (he1-12) -- (w1);
	
	\path[line width=1.5pt]

	(v1) edge [edge, uhhdarkgrey, bend right] (v0)
	(v2) edge [edge, uhhdarkgrey, bend right] (v1)
	(v3) edge [edge, uhhdarkgrey, bend right] (v2)
	(v3) edge [edge, uhhdarkgrey, loop right, looseness=8] (v3)

	(v2) edge [edge, grey, bend left] (v0)	
	(v3) edge [edge, grey, bend left] (v1)	
	(v3) edge [edge, grey, bend left] (v2)
	(v3) edge [edge, grey, loop right, looseness=12] (v3)
	;
 
	\begin{pgfonlayer}{background}
		\draw[rounded corners, fill=grey!20] (-1,.75) -- (2.3,.75) -- (2.3,-.75) -- (-2,-.75) -- (-2,.75) -- (-1,.75);
	\end{pgfonlayer}

\end{tikzpicture}%
			}%
			\caption{An augmented hypernetwork with a core consisting of a classical feed-forward network with four nodes}
			\label{fig:4th_ex}
		\end{subfigure}
	\end{minipage}
	\hfill
	\begin{minipage}{.6\linewidth}
		\centering
		\begin{subfigure}{\linewidth}
			\includegraphics[width=\linewidth]{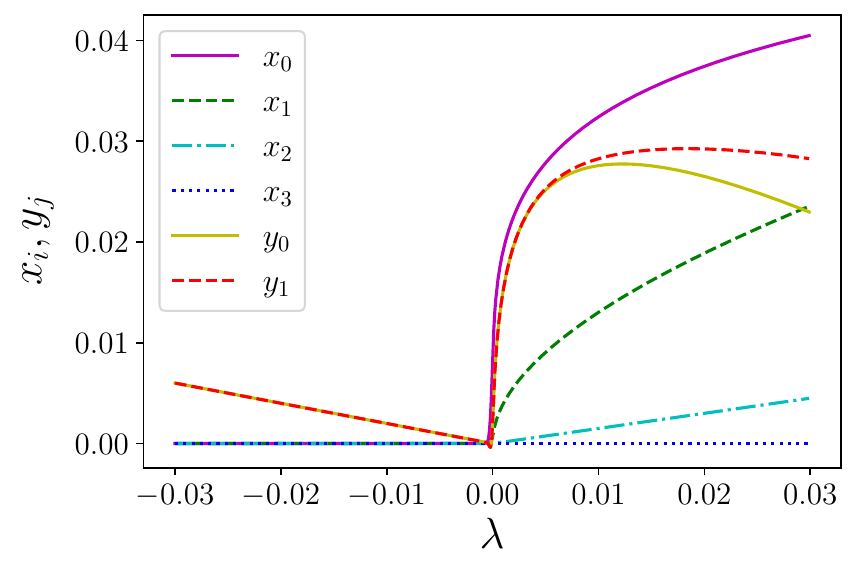}
			\caption{The stable branches of a synchrony-breaking bifurcation.}
			\label{fourthnum-a}
		\end{subfigure} \\
		\begin{subfigure}{\linewidth}
			\includegraphics[width=\linewidth]{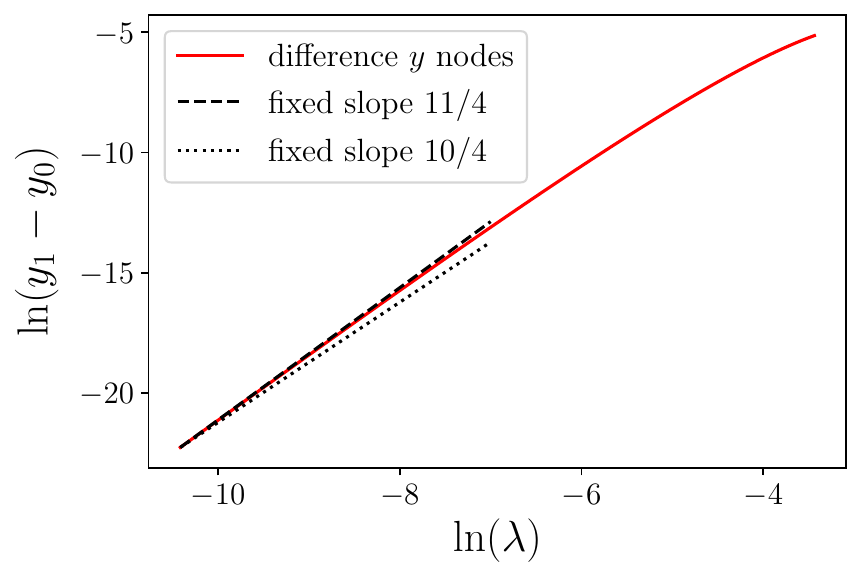}
			\caption{A log–log plot of the difference between the $y$-nodes.}
			\label{fourthnum-b}
		\end{subfigure}
	\end{minipage}    
	\caption{An augmented hypernetwork with a core consisting of a classical feed-forward network with four nodes, shown within the grey box. We have left out self-loops corresponding to self-influence of each node. In contrast to the previous examples, for this hypernetwork we have $k+1=4$, so that there are $(k+1)!=24$ hyperedges of order $k=3$. Also shown  are numerically obtained bifurcation diagrams. The dashed black line segment in the log-log plot has fixed slope $\frac{11}{4}$, whereas the dotted segment has fixed slope $\frac{10}{4}$.}
	\label{fourthnum}
\end{figure}  

\begin{ex}
	Finally, we consider the augmented hypernetwork of \Cref{fig:4th_ex}, which has as its core a classical feed-forward network with four nodes. 
	More precisely, the nodes in the core evolve according to the ODEs
	\begin{equation}
		\label{system4feedfor}
		\begin{array}{ll}
			\dot x_0 & = G(x_0, \textcolor{uhhdarkgrey}{x_1}, \textcolor{grey}{x_2}; \lambda)\, , \\
			\dot x_1 & = G(x_1, \textcolor{uhhdarkgrey}{x_2}, \textcolor{grey}{x_3}; \lambda) \, , \\
			\dot x_2 & = G(x_2, \textcolor{uhhdarkgrey}{x_3}, \textcolor{grey}{x_3}; \lambda)\, , \\
			\dot x_3 & = G(x_3, \textcolor{uhhdarkgrey}{x_3}, \textcolor{grey}{x_3}; \lambda)\, . \\
		\end{array}
	\end{equation}
	It is known that this system generically supports steady-state bifurcations in which stability passes from a fully synchronous branch to one in which
	\begin{align*}
		x_0(\la) \sim \la^{1/4}\, , \,\,
		x_1(\la) \sim \la^{1/2} \text{ and } \,
		x_2(\la), x_3(\la), x_2(\la)- x_3(\la) \sim \la\, ,
	\end{align*}
	see \cite{vondergracht2022} and \cite{rink2013amplified}.  This unusually fast rate of synchrony breaking is also referred to as \emph{amplification}.
	It follows that 
	\begin{equation*}
		\bar{p} = \frac{1}{4} + \frac{1}{4} + \frac{1}{4} + \frac{1}{2} + \frac{1}{2} + 1 = \frac{11}{4}\, ,
	\end{equation*}
	so that \Cref{mainbifurcation} predicts a reluctant steady state branch with 
	\begin{equation*}
		y_0(\lambda) - y_1(\lambda) \sim \lambda^{\frac{11}{4}}\, .
	\end{equation*}
	This unusual growth rate is verified numerically in \Cref{fourthnum}, which  was obtained by numerically integrating the augmented system for
	\begin{equation*}
		G(X_0, X_1, X_2; \lambda) =  10X_1 -20X_2 + 15\la X_0 -100X_0^2\, 
	\end{equation*}
	and
	\begin{align*}
		F\left(Y, \bigoplus_{\sigma \in S^0_{4}}  {\bf X}_{\sigma};\la\right) &= -0.01\sum_{\sigma \in S^0_{4}}
		h(120 X_{\sigma,1}+40X_{\sigma, 2}-100X_{\sigma,3})
		-5Y-\la \, ,
	\end{align*}
	where ${\bf{X}}_{\sigma} = (X_{\sigma,1}, \dots, X_{\sigma,4}) \in \R^4$ and with $h$ given by  \Cref{num:res-h}. \\ 
	\indent \Cref{fourthnum-a} shows the components of the stable branches, whereas \Cref{fourthnum-b} is a log-log plot of the difference of the $y$-components, for positive values of $\la$. In this latter picture, the dashed black line segment has fixed slope $11/4$, indicating that indeed $y_0(\lambda) - y_1(\lambda) \sim\lambda^{\frac{11}{4}}$. 
	For comparison we also plotted the dotted black line segment, which has fixed slope $10/4$, and which does not fit as well for low values of $\lambda$. Both figures were made by forward integrating the vector field for the augmented system from the point $$(x_0, \dots, x_3, y_0, y_1) = (-0.001, -0.002, -0.003, -0.004, 0.001, 0.002)$$ in phase space, for various values of $\lambda$ and with time steps of $0.1$. 
	For \Cref{fourthnum-a} this was done up to $t=2000$ and with $600$ equidistant values of $\la \in [-0.03, 0.03]$. 
	For \Cref{fourthnum-b} this was up to $t=20000$ and for $100$ equidistant values of $\ln(\la) \in [\ln(0.00003), \ln(0.03)]$.
\end{ex}

\begin{remk}
	\label{remk:secondclosing}
	\Cref{example:largefund} shows that even classical (dyadic) networks may generically support  reluctant synchrony breaking bifurcations. 
	However, we are not aware of any method to design networks that, for instance, break synchrony up to some prescribed degree in $\lambda$.  
	For hypernetworks, the augmented hypernetwork construction makes this design problem  more tractable. In fact, we show below that one may construct hypernetworks which support generic reluctant synchrony breaking to arbitrarily high order.  
	
	To illustrate how one can create hypernetworks with an arbitarily high order of reluctant synchrony breaking, we return to the five-node augmented hypernetwork of \Cref{vectorfieldex1}. 
	In \Cref{ex:leadingreluctant}, we observed that the three-node core of this hypernetwork supports a steady state branch in which $x_i(\la) = D_i\la + \mathcal{O}(|\lambda|^2)$ for some mutually distinct $D_i \in \R$. 
	Using the same notation as in  \Cref{ex:leadingreluctant}, we   expand the response function $F$ for the $y$-nodes as 
	\begin{equation}\label{expressionsforFnnfG--2}
		\begin{split}
			& F(Y, (X_0, X_1), (X_2, X_3), (X_4, X_5); \lambda )  \\   
			& = aY + bX_0 + cX_1 + bX_2 + cX_3 + bX_4 + cX_5 + d\lambda \\   
			& \phantom{=} + \mathcal{O}(\|(Y,X_0, \dots, X_5;\lambda)\|^2) \, ,
		\end{split}
	\end{equation}
	and we recall that we found a reluctant steady state branch in the augmented hypernetwork with asymptotics 
	\begin{equation}\label{y0non-synns--2}
		y_0(\lambda) = \frac{-(b+c)(D_0 + D_1 + D_2) - d}{a}\lambda  + \mathcal{O}(|\lambda|^2)\, 
	\end{equation}
	and
	\begin{equation}\label{y1non-synns--2}
		y_1(\lambda) = \frac{-(b+c)(D_0 + D_1 + D_2) - d}{a}\lambda  + \mathcal{O}(|\lambda|^2)\, .
	\end{equation}
	To this augmented hypernetwork we can now add another node of the same type as the $y$-nodes, with corresponding variable $y_2$. 
	We couple it to the nodes in the core in such a way that
	\begin{equation}\label{y2eeq--2}
		\dot{y_2} = F(y_2, (x_0, x_1), (x_1, x_0), (x_0, x_0); \lambda)\, .
	\end{equation}
	The aforementioned branch of steady states is then supported by this larger hypernetwork as well, where in addition
	\begin{equation}\label{y2non-synns--2}
		y_2(\lambda) = \frac{-(b+c)(D_0 + D_1 + D_0) - d}{a}\lambda  + \mathcal{O}(|\lambda|^2)\, ,
	\end{equation}
	as can be seen using \Cref{expressionsforFnnfG--2}.
	To summarise, we now have a branch where the three $y$-nodes satisfy
	\begin{align}
		y_0(\la) &= E_0\la + \mathcal{O}(|\la|^2) \, , \, \, y_1(\la) = E_0\la + \mathcal{O}(|\la|^2) \, ,\\ \nonumber
		y_2(\la) &= E_2\la + \mathcal{O}(|\la|^2) \text{ and } \, y_0(\la) - y_1(\la) \sim \la^p
	\end{align}
	for some $p>1$ (in this particular case $p=3$). 
	Moreover, from Equations \eqref{y0non-synns--2}, \eqref{y1non-synns--2} and \eqref{y2non-synns--2} we see that $E_0 \not= E_2$, as $D_0 \not= D_2$ by assumption. \\
	\indent We may now use the three $y$-nodes as the core for another augmented hypernetwork, say by adding two $z$-nodes of a new type (c.f. \Cref{remk:firstclosing}). 
	We also add a third $z$-node, precisely as we did with the third $y$-node. 
	That is, we set
	\begin{equation}\label{z2eeq--2}
		\begin{split}
			\dot{z_0} &= \tilde{F}(z_0, (y_0, y_1), (y_1, y_2), (y_2, y_0); \lambda)\, ,\\
			\dot{z_1} &= \tilde{F}(z_1, (y_0, y_2), (y_1, y_0), (y_2, y_1); \lambda)\, ,\\
			\dot{z_2} &= \tilde{F}(z_2, (y_0, y_1), (y_1, y_0), (y_0, y_0); \lambda)\, ,
		\end{split}
	\end{equation}
	where $\tilde{F}$ is a response function for the $z$-nodes. Just as before, we will then find 
	\begin{align}
		z_0(\la) &= E'_0\la + \mathcal{O}(|\la|^2) \, , \, \, \, z_1(\la) = E'_0\la + \mathcal{O}(|\la|^2) \, ,\\ \nonumber
		z_2(\la) &= E'_2\la + \mathcal{O}(|\la|^2) \text{ and } z_0(\la) - z_1(\la) \sim \la^{p+2}
	\end{align}
	for some generically nonzero $E'_0, E'_2 \in \R$. The new power ${p+2}$ follows from \Cref{mainbifurcation}, as
	\begin{equation}
		y_2(\la) - y_0(\la) \sim \la\, , 
		\quad y_2(\la) - y_1(\la) \sim \la \, , \quad
		y_0(\la) - y_1(\la) \sim \la^p\, ,
	\end{equation}
	which holds because $E_0 \not= E_2$. We may also argue that $E'_0 \not= E'_2$ in precisely the same way that we argued that $E_0 \not= E_2$. \\ 
	\indent 
	This shows that by iteratively growing the augmented hypernetwork, we may increase the order (in $\lambda$) of reluctancy of the reluctant steady state branch. In other words, we may design hypernetworks with an arbitrarily high order of reluctant synchrony breaking. Concretely, our example shows that we can arrange for $p = 3, 5, 7, \dotsc$. 
	It is also clear from this construction that the resulting reluctant branch may be assumed stable.
\end{remk}

\section{Discussion}
\label{sec:discussion}
In \cite{vonderGracht.2023}, the authors introduced a mathematical framework to capture  higher order interactions in network dynamical systems, thereby generalising the analogous set-up for classical (dyadic) networks \cite{Golubitsky.2004, Golubitsky.2006, Field.2004}.
It is observed in \cite{vonderGracht.2023} that, unlike for dyadic networks,  synchronisation in these hypernetwork dynamical systems is governed by higher order (nonlinear) terms in the equations of motion. This suggests that hypernetwork systems may display interesting phenomena that cannot be observed in dyadic networks.   
In particular, the authors of \cite{vonderGracht.2023} construct a hypernetwork system that shows numerical evidence of a bifurcation scenario in which two nodes break synchrony at an unusually high order in the bifurcation parameter.

In this paper we provide a general method to  construct  hypernetwork systems  displaying such ``reluctant'' synchrony breaking, and we give a rigorous mathematical proof that these bifurcations  occur generically in these systems. We also give an analytical expression for the order in the bifurcation parameter at which the synchrony breaking occurs. 

Even though reluctant synchrony breaking is not impossible for dyadic networks, there is currently no understanding of what causes their (rare) occurrence. 
A method for constructing networks to achieve reluctant synchrony breaking is likewise lacking. 
In this paper, we show that this design problem is much more tractable for hypernetwork systems, and that the phenomenon appears to be significantly more common when higher order interactions are present.
This sheds new light on the role that such higher order interactions play in various natural systems, and on their potential for applications in engineering.

In particular, we see interesting parallels with the concept of \emph{homeostasis}. 
This term refers to the ability of living organisms to keep their internal states approximately stable when external conditions are changed. 
A well-known example is the ability of warm-blooded animals to regulate their body temperature across a wide range of environmental temperatures \cite{Golubitsky.2017b, wang2021structure}.
For dynamical systems with a distinguished input parameter $\mathcal{I}$ and output function $\mathcal{X}$, an (infinitesimal) homeostasis point is an input parameter value $\mathcal{I} = \mathcal{I}_0$  at which $\frac{d\mathcal{X}}{d\mathcal{I}} = 0$, with  higher order derivatives  possibly  vanishing too \cite{Golubitsky.2023}.
In our set-up, the output function is the difference between the states of two nodes at equilibrium, while the input parameter is a bifurcation parameter $\lambda$. Thus, reluctant synchrony breaking can be interpreted as a form of ``synchrony homeostasis''. 
This phenomenon may well occur within living systems, as synchrony has been linked for instance to gene-regulatory processes \cite{morone2020fibration} and memory formation in the brain \cite{jutras2010synchronous}. 
Of course, such claims would have to be verified in the relevant biological systems. Nevertheless, it is striking that higher order interactions can slow down de-synchronisation, perhaps to make a system more resilient to variations in external conditions. 

\section{Acknowledgements}
We thank Martin Golubitsky and Ian Stewart for enlightening discussions,  as well as Tiago Pereira and Deniz Eroglu for helpful remarks on a draft of this manuscript. 

S.v.d.G. was partially funded by the Deutsche Forschungsgemeinschaft (DFG, German Research Foundation)--–453112019. E.N. was partially supported by the Serrapilheira Institute (Grant No. Serra-1709-16124). B.R. gratefully acknowledges the hospitality and financial support of the Sydney Mathematical Research Institute. 

\section{Data availability}
The code used to generate the numerical results can be found in \cite{githubcodereluctantptII}.

\addcontentsline{toc}{section}{References}
\printbibliography

\end{document}